\documentclass[psamsfont]{amsart}

\usepackage{graphicx}
\usepackage{amssymb}
\usepackage{xypic}
\usepackage{hyperref}
\usepackage{psfrag}
\usepackage{amsmath,amssymb}
\usepackage{dsfont}\let\mathbb\mathds

\newtheorem{theo}{Theorem}[section]
\newtheorem{coro}[theo]{Corollary}
\newtheorem{lemm}[theo]{Lemma}
\newtheorem{prop}[theo]{Proposition}
\theoremstyle{definition}
\newtheorem{defi}[theo]{Definition}
\newtheorem{rema}[theo]{Remark}

\def\del{\partial}              
\def\eps{\varepsilon}           

\def\bC{\mathbb{C}}         
\def\bR{\mathbb{R}}         
\def\bZ{\mathbb{Z}}         

\def\mA{\mathcal{A}}            
\def\mC{\mathcal{C}}            
\def\mI{\mathcal{I}}            
\def\mL{\mathcal{L}}            
\def\mM{\mathcal{M}}            
\def\mP{\mathcal{P}}            
\def\mQ{\mathcal{Q}}            
\def\mR{\mathcal{R}}            
\def\mS{\mathcal{S}}            
\def\mW{\mathcal{W}}            

\def\Ham{\mathrm{Ham}(M,\omega)}            

\def\ra{\rightarrow}
\def\lra{\longrightarrow}
\def\co{\colon\thinspace}

\def\skc#1#2{\underline{\sigma}^{#1}_{X,[l]}(#2;L,H,J,\eta)}          
\def\sic#1#2{\overline{\sigma}^{#1}_{X,[l]}(#2;L,H,J,\eta)}           
\def\kc#1#2{\underline{c}^{#1}_{X,[l]}(#2;L,L')}                      
\def\ic#1#2{\overline{c}^{#1}_{X,[l]}(#2;L,L')}                       

\def\hop{{\noindent}}

\begin{document}

\title[Lagrangian spectral invariants]
      {Spectral invariants in Lagrangian Floer theory}
\author[R\'emi Leclercq]{R\'emi Leclercq}
\address{R.L.: Universit\'e de Montr\'eal\\
Department of Mathematics and Statistics\\
CP 6128 Succ. Centre Ville\\
Montr\'eal, QC H3C 3J7\\
Canada}
\email{leclercq@dms.umontreal.ca}
\subjclass[2000]{Primary 57R17; Secondary 53D12 53D40 55T10}
\keywords{Spectral invariants, spectral sequence, Lagrangian submanifolds, Lagrangian Floer homology, Hofer's distance.}

\maketitle

\begin{abstract}
Let $(M,\omega)$ be a symplectic manifold compact or convex at infinity. Consider a closed Lagrangian submanifold $L$ such that $\omega |_{\pi_2(M,L)}=0$ and $\mu|_{\pi_2(M,L)}=0$, where $\mu$ is the Maslov index. Given any Lagrangian submanifold $L'$, Hamiltonian isotopic to $L$, we define Lagrangian spectral invariants associated to the non zero homology classes of $L$, depending on $L$ and $L'$. We show that they naturally generalize the Hamiltonian spectral invariants introduced by Oh and Schwarz, and that they are the homological counterparts of higher order invariants, which we also introduce here, via spectral sequence machinery introduced by Barraud and Cornea. These higher order invariants are new even in the Hamiltonian case and carry strictly more information than the classical ones. We provide a way to distinguish them one from another and estimate their difference in terms of a geometric quantity.
\end{abstract}



\section{Introduction}

Hamiltonian spectral invariants for Floer homology have been defined by Schwarz \cite{Schwarz00} and Oh \cite{Oh05}, following the work of Viterbo~\cite{Viterbo92}. Schwarz worked in the symplectically aspherical case and Oh has treated the general case. Oh \cite{Oh97} and Milinkovi\'c \cite{Milinkovic00} also defined a Lagrangian version of these invariants for cotangent bundles.

In this paper, we define a new type of Lagrangian spectral invariants. The novelty is twofold. Firstly, we use the PSS morphism to define them (this is the same approach as Schwarz's in the Hamiltonian case). Secondly, our invariants are of higher order, in the sense that they are associated to each non zero class $\alpha\in E^{r}_{p,q}(L;X)$ where $E(L;X)$ (also denoted $E_{X}(L)$) is a purely topological Serre spectral sequence associated to $L$ (and an auxiliary space $X$) via a construction due to Barraud and Cornea~\cite{BarraudCornea07}. These invariants are new even in the Hamiltonian case, where they extend the classical invariants of Oh and Schwarz which are recovered for $r=2$. Moreover they carry strictly more information than the classical invariants, even in the Morse case. We study some properties of these invariants and we describe, in particular, some geometric way to distinguish them -- this part is again new even when applied to the usual Hamiltonian case. Here are now our main results described in more detail.\\

Let $(M,\omega)$ be a symplectic manifold, compact or convex at infinity, and $L$ be a closed Lagrangian submanifold such that $\omega|_{\pi_2(M,L)}=0$ and $\mu|_{\pi_2(M,L)}=0$ with $\mu$ the Maslov index. Let $(H,J)$ be a regular pair, formed of a Hamiltonian function and an almost complex structure. Let $\phi_H^1$ be the time--$1$ diffeomorphism induced by $H$. We denote by $L':=\phi_H^1(L)$ and we assume $L\cap L'$ transverse. We will start with the definition of order $2$ (or homological) spectral invariants since this construction will be more familiar to the reader. All homology theories used here are with coefficients in $\bZ_2$.

In our setting, the Lagrangian Floer complex $CF_*(L,L';H,J)$ is well-defined and admits a natural filtration by action, $CF^\nu_*(L,L';H,J)$. In the Hamiltonian case, Piunikhin, Salamon and Schwarz~\cite{PiunikhinSalamonSchwarz96} introduced a particular morphism between the Morse and Floer complexes inducing an isomorphism in homology. It has been adapted to the Lagrangian setting by Kati\'c and Milinkovi\'c~\cite{KaticMilinkovic05} for cotangent bundles and in more generality by Barraud and Cornea~\cite{BarraudCornea05}, and Albers~\cite{Albers05}. We denote it by $\phi^H_f\co HM_*(L;f,g)\ra HF_*(L,L;H,J)$, with $HM_*(L;f,g)$ the Morse homology of $L$ induced by a Morse--Smale pair $(f,g)$. Let $\eta$ be a path from $L$ to $L'$.
\begin{defi}\label{defi:si}
Let $\alpha$ be a non zero (Morse) homology class of $L$. Its associated relative homological Lagrangian spectral number is given by
\begin{align*}
\sigma_L(\alpha;H,J,\eta;f,g):= \inf \{ \nu \in\bR |\, \phi_f^H(\alpha) \in {\rm Im}(i^\nu_*) \}
\end{align*}
where $i^\nu_*=H_*(i^\nu)$ is induced by the inclusion $CF^\nu_*(L,L';H,J)\hookrightarrow CF_*(L,L';H,J)$. Its (absolute) {\rm homological Lagrangian spectral number} is defined as
\begin{align*}
c_L(\alpha;H,J;f,g):= \sigma_L(\alpha;H,J,\eta;f,g) - \sigma_L(1;H,J,\eta;f,g)
\end{align*}
with $1$ denoting the generator of $H_0(L)$.
\end{defi}
Along these lines, homological Lagrangian spectral numbers will refer to this absolute version. The distinction {\it relative} vs {\it absolute} reflects the fact that in the Lagrangian version of Floer theory, the action functional is defined up to translation by a constant due to the choice of $\eta$, in some fixed homotopy class of paths from $L$ to $L'$.

The construction of these numbers requires various choices but, as we shall see, they are strongly invariant. Recall that Hofer's distance between two Hamiltonian isotopic Lagrangian submanifolds, denoted $\nabla(L,L')$, is the infimum of the energies of Hamiltonian diffeomorphisms which carry $L$ to $L'$.
\begin{theo}\label{theo:independ}
The homological Lagrangian spectral numbers depend only on $L$ and $(\phi^1_H)^{-1}(L)$. In other words, for any Lagrangian submanifold $L'$, Hamiltonian isotopic to $L$, and any non zero homology class $\alpha\in H_*(L)$, we may define, its associated {\rm Lagrangian spectral invariant}
\begin{align*}
c(\alpha;L,L'):= c_L(\alpha;H,J;f,g)
\end{align*}
with $H$ such that $\phi^1_H(L')=L$.

Moreover, the application $c(\alpha;L,-)$ is continuous with respect to Hofer's distance on the set of all Lagrangian submanifolds Hamiltonian isotopic to $L$.
\end{theo}
\hop The most delicate part of the proof of invariance consists in showing the commutativity of a certain diagram (\ref{diag}). This is shown via methods inspired by Seidel~\cite{Seidel97} and which have been used in the Hamiltonian case by McDuff and Salamon~\cite{McDuffSalamon04}. \\

These quantities are a natural generalization of the Hamiltonian spectral invariants. Indeed, let $\Delta$ be the diagonal in $M\times M$, and $\underline{\alpha}\in H_*(\Delta)$ be the homology class corresponding to $\alpha\in H_*(M)$, via the obvious isomorphism. For $\phi\in \Ham$, we denote by $\rho(-;\phi)$ its associated spectral invariant and by $\Gamma_{\phi}$ its graph. We have
\begin{align}\label{prop:generalisation}
c(\underline{\alpha};\Delta,\Gamma_{\phi})= \rho(\alpha;\phi) - \rho(1;\phi)
\end{align}
where $1$ denotes the generator of $H_0(M)$. This is shown via an isomorphism due to Biran, Polterovich and Salamon \cite{BiranPolterovichSalamon03}. The Lagrangian spectral invariants, as defined before, are easily seen to satisfy two properties, similar to the ones valid in the Hamiltonian case, and listed below.
\begin{prop}\label{prop:dual}
Let $\alpha\in H_k(L)$, $\alpha\neq 0$. We denote by $\alpha'\in H_{n-k}(L)$ the $Hom$--dual of its Poincar\'e dual class.
\begin{itemize}
\item[{\it (1.)}] For $\phi\in \Ham$, $c(\alpha;L,\phi(L))= c([L];L,\phi^{-1}(L))-c(\alpha';L,\phi^{-1}(L))$
\item[{\it (2.)}] $0\leq c(\alpha;L,L')\leq \nabla(L,L')$, with strict inequalities as soon as $0<k<n$
\end{itemize}
\end{prop}
\hop Viewing the Hamiltonian spectral invariants as particular Lagrangian ones via (\ref{prop:generalisation}), assertion {\it (2.)} improves the well-known upper bound given by Hofer's norm of $\phi$. We will illustrate this fact via a family of Hamiltonian diffeomorphisms given by Ostrover~\cite{Ostrover03}.\\

The Lagrangian spectral invariants defined above are actually the simplest particular cases of higher order spectral invariants which we now discuss. These higher order invariants are denoted $\ic{r}{\alpha}$ and $\kc{r}{\alpha}$ where $l\co L\ra X$ is a map with the space $X$ simply-connected, $[l]$ is its homotopy class and $\alpha$ is any non  zero element of the $r$--th page of a spectral sequence $E_X(L)$ defined as follows. We consider the path-loop fibration $\Omega X \to PX \to X$. We pull-back this fibration over the map $l$ and we let $E_{X}(L)$ be the Serre spectral sequence of this pull-back fibration.

It has been shown by Barraud and Cornea \cite{BarraudCornea05} that this spectral sequence can be related by a PSS type isomorphism to a spectral sequence induced by a natural filtration of an enriched Floer type complex. In view of this, Definition \ref{defi:si} can be adapted to this element $\alpha$ (see definitions \ref{defi:si2} and \ref{defi:si2bis}).

The second page of $E_X(L)$ is purely homological, in the sense that $E^2_{p,q}(L;X)\simeq H_q(\Omega X)\otimes H_p(L)$. Thus, any element $\alpha\neq 0$ in $H_p(L)$ can be identified with the element $1\otimes\alpha$ in $E^2_{p,0}(L;X)$ where $1$ denotes the generator of $H_0(\Omega X)$. For such an element, all three Lagrangian invariants coincide:
\begin{align*}
\ic{2}{1\otimes\alpha}=\kc{2}{1\otimes\alpha}= c(\alpha;L,L').
\end{align*}
The higher order invariants generalize the classical spectral invariants (in the previous sense) and carry strictly more information (even in the Morse case). We will illustrate this fact with explicit computations in a particular case. As suggested by our notation, given $X$, $[l]$ and $r$, $\ic{r}{-}$ and $\kc{r}{-}$ have the same invariance property as their homological counterparts.

The main property that will be of interest to us here is that $E_X(L)$ provides a way to distinguish among these invariants, in terms of a quantity related to the geometry of $L$ and $L'$. This geometric constant is defined, for every pair of transverse Lagrangian submanifolds, as follows. Let $x$ be an element of $L\cap L'$. There exist a real number $\eps >0$ and an embedding, $e^x_\eps$, of the ball $B(0,\eps)\subset\bC^n$ in $M$ such that
\begin{equation}\label{boules}
\begin{split}
{\rm i.} \;\; & (e^x_\eps)^*(\omega) = \omega_0 \;\;\mbox{ and }\;\; e^x_\eps(0)=x, \\
{\rm ii.} \;\; & (e^x_\eps)^{-1}(L) = \bR^n\cap B(0,\eps)\;\; \mbox{ and } \;\; (e^x_\eps)^{-1}(L') = i\bR^n\cap B(0,\eps).
\end{split}
\end{equation}
This gives rise to the following definition.
\begin{defi}\label{def:rll}
Given $L$ and $L'$ two compact, transverse Lagrangian submanifolds, we define
\begin{align*}
r(L,L') := {\rm sup} \left\{ \eps >0 \, \left| \;
\begin{array}{l}
\forall x\in L\cap L',\; \exists\, e^x_\eps \mbox{ satisfying conditions {\rm (\ref{boules})}}, \\
\mbox{such that }\; x\neq y \;\Rightarrow\; {\rm Im}\, e^x_\eps \cap {\rm Im}\, e^y_\eps = \emptyset
\end{array}
\right. \right\}.
\end{align*}
\end{defi}
This definition is inspired by the geometric distance introduced by Barraud and Cornea \cite{BarraudCornea07}. Notice that, since $L\cap L'$ is a finite set, $r(L,L')>0$. With these conventions we have:
\begin{theo}\label{theo:main}
If $d^r(\alpha)=\beta\neq 0$ in $E_X(L)$, then
\begin{align*}
\kc{r}{\alpha} - \kc{r}{\beta} &\geq 0 \;\;{\rm and} \\
\ic{r}{\alpha} - \kc{r}{\beta} &\geq \frac{\pi r(L,L')^2}{2}
\end{align*}
for any Lagrangian submanifold $L'$, Hamiltonian isotopic (and transverse) to $L$.
\end{theo}

This property of Lagrangian spectral invariants of higher order has the following immediate consequence for the invariants of Definition \ref{defi:si}.
\begin{coro}\label{theo:ss}
Let $X$ be $(r-1)$--connected and let $\{x_i\}$ be a basis of $H_{r-1}(\Omega X)$. If there exists a non trivial differential $d^r\co H_0(\Omega X)\otimes H_p(L)\simeq E^r_{p,0}\ra E^r_{p-r,r-1}$ of $E_X(L)$, namely $d^r(1\otimes\alpha)= \sum x_i\otimes \beta_i \neq 0$, then for all Lagrangian submanifold $L'$, Hamiltonian isotopic to $L$,
\begin{align*}
\forall i,\; c(\alpha;L,L')-c(\beta_i;L,L') \geq \frac{\pi r(L,L')^2}{2}.
\end{align*}
\end{coro}
\hop This corollary has interesting consequences.
\begin{coro}\label{prop:prop}
Let $\cdot$ denote the intersection product and let $[L]$ be the fundamental class of $L$.
\begin{itemize}
\item[{\it (1.)}] For $\alpha\in H_k(L)$ and $\beta\in H_*(L)$, with $1<k<n-1$ and $\alpha\cdot\beta\neq 0$, we have
\begin{align*}
c(\alpha\cdot\beta; L,L') \leq c(\beta; L,L') - \frac{\pi r(L,L')^2}{2}
\end{align*}
\item[{\it (2.)}] If $\alpha\neq 0$ is an element of $H_k(L)$ with $1<k<n-1$, then it holds:
\begin{align*}
\frac{\pi r(L,L')^2}{2} \leq c(\alpha; L,L') \leq c([L]; L,L') - \frac{\pi r(L,L')^2}{2}
\end{align*}
\end{itemize}
\end{coro}
As we shall see, in this purely homological case, there is another way to prove these properties by adapting to the Lagrangian setting the proofs valid in the Hamiltonian case (and used by Schwarz \cite{Schwarz00}) and combining this with the definition of $r(L,L')$. The key point of this alternate proof is an adequate description of a certain module structure on $HF_*(L,L')$ over $H_*(L)$ endowed with its intersection product. The cases $k=1$ and $k=n-1$ of the previous corollary are proved here via this method, in the sense that assertions {\it (1.)} and {\it (2.)} hold as long as $0<k<n$. In particular, for all non zero $\alpha\in H_k(L)$, $0 \leq c(\alpha; L,L') \leq c([L]; L,L')$.

Finally we get, as an obvious consequence of this extension of Corollary \ref{prop:prop} {\it (1.)} and Proposition \ref{prop:dual} {\it (2.)}, a lower bound for Hofer's distance between two transverse, Hamiltonian isotopic Lagrangian submanifolds in terms of our geometric constant and the cup-length of $L$.
\begin{coro}\label{coro:cup-length}
Let ${\rm cl}(L)$ be the cup-length of $L$ (with $\bZ_2$ coefficients), we have
\begin{align*}
\nabla(L,L') \geq {\rm cl}(L)\cdot\frac{\pi r(L,L')^2}{2}
\end{align*}
for all Lagrangian submanifold $L'$, Hamiltonian isotopic and transverse to $L$.
\end{coro}

\hop{\bf Organization of the paper.} In Section \ref{section1}, we first quickly recall the construction of Morse and Floer homologies and three kinds of comparison morphisms, namely the classical comparison morphism for Lagrangian Floer homology, a Lagrangian version of the PSS morphism and the naturality morphism.

In Section \ref{sec:demo}, we construct our Lagrangian spectral numbers and we prove Theorem \ref{theo:independ}. The invariance part relies on two distinct steps. First we show that a certain diagram commutes. In order to do this, we prove that the Lagrangian Floer homology $HF_*(L,L')$ can be viewed as a module over the homology of $L$, endowed with its intersection product. This algebraic structure is an adaptation of a construction made by Floer \cite{Floer89} in the Hamiltonian case. Then we show that the naturality and Lagrangian PSS morphisms preserve this structure. The desired commutativity follows immediately. In a second step, we adapt methods used by Schwarz in the Hamiltonian case, to get an estimate on Lagrangian spectral numbers. This estimate ends the proof of their invariance and implies the continuity part of Theorem \ref{theo:independ}. Finally we show how the Hamiltonian spectral invariants can be viewed as particular Lagrangian ones, by proving (\ref{prop:generalisation}).

In Section \ref{sec:3}, we first sketch the construction of the Barraud--Cornea spectral sequence and its relation with the Serre spectral sequence of the path-loop fibration of $L$. Then we introduce the Lagrangian spectral invariants of higher order and prove Theorem \ref{theo:main}. We also produce explicit computations of the higher order spectral invariants and show, via this example, that they carry more information than the classical invariants.

Finally (in Section \ref{sec:4}), we come back to the homological Lagrangian spectral invariants. Corollary \ref{theo:ss} follows immediately from the results obtained on the higher order invariants. We show how it implies Corollary \ref{prop:prop}. Then, we prove Proposition \ref{prop:dual} and Corollary \ref{coro:cup-length}.\\

\hop{\bf Acknowledgements.} The author would like to deeply thank his supervisor, Octav Cornea, for hours of extremely fruitful and enlightening discussions. He would also like to thank Matthias Schwarz for valuable conversations, Peter Albers for numerous interesting comments and discussions and Leonid Polterovich for useful questions.


\section{Recalls and notation}\label{section1}

\subsection{Morse and Floer homologies}
\subsubsection{Morse homology}

We briefly sketch this well-known construction. We work over $\bZ_2$, and fix a closed manifold $L^n$, a metric $g$ and a Morse function $f$ such that the pair $(f,g)$ is Morse--Smale. We denote the gradient of $f$ with respect to $g$ by $\nabla f$. The gradient flow of $-f$, defined by $\frac{d}{dt} \gamma_t +\nabla f(\gamma_t)=0$, induces stable and unstable manifolds for every critical point $p$:
\begin{align*}
\mW^s_p(f,g) &:=\{ x\in M \,|\; \lim_{t\ra +\infty} \gamma_t(x)=p \} \mbox{ and} \\
\mW^u_p(f,g) &:=\{ x\in M \,|\; \lim_{t\ra -\infty} \gamma_t(x)=p \}.
\end{align*}

Let $i_f(p)$ be the Morse index of $p$. As $f$ is Morse--Smale, the connecting manifold of two critical points $p$ and $q$, $\widehat{\mM}_{p,q}(f,g):= \mW^u_p(f,g)\cap \mW^s_q(f,g)$ is a smooth ($i_f(p)-i_f(q)$)--dimensional manifold. As $\bR$ acts on this space, one can define $\mM_{p,q}(f,g):=\widehat{\mM}_{p,q}(f,g)/\bR$. In general, these moduli spaces are not compact but admit a compactification such that
\begin{align}\label{compact}
\del \overline{\mM}_{p,q}(f,g)=\bigcup_{r\in {\rm Crit}(f)} \overline{\mM}_{p,r}(f,g)\times \overline{\mM}_{r,q}(f,g).
\end{align}

Let ${\rm Crit}_k(f)$ be the set of critical points of $f$ with Morse index $k$. We define free $\bZ_2$--modules $CM_k(L;f,g):=\langle {\rm Crit}_k(f)\rangle_{\bZ_2}$ and a differential $\del$ by the formula
\begin{align*}
\del p=\sum_{q\in {\rm Crit}_{k-1}(f)} \#_2 \mM_{p,q}(f,g)\cdot q
\end{align*}
where $\#_2 \mM_{p,q}(f,g)$ denotes the mod 2 cardinal of $\mM_{p,q}(f,g)$. The compactification (\ref{compact}) implies that $\del^2=0$. Thus $(CM_*(L;f,g),\del)$ is a chain complex whose homology
\begin{align*}
HM_*(L):=H(CM_*(L;f,g),\del)
\end{align*}
is the Morse homology of $L$. It is a well-known fact that it does not depend on the Morse--Smale pair $(f,g)$ and is the (cellular) homology of $L$.

\subsubsection{Lagrangian Floer homology}

This theory is exposed in Floer~\cite{Floer88}, Oh~\cite{Oh93}. Let $(M,\omega)$ be a symplectic manifold, compact or convex at infinity. Let $L$ and $L'$ be two closed, Hamiltonian isotopic Lagrangian submanifolds of $M$ such that
\begin{align}\label{context}
\omega|_{\pi_2(M,L)}=0\;\;{\rm and}\;\; \mu|_{\pi_2(M,L)}=0
\end{align}
where $\mu$ is the Maslov index (its definition is recalled below). We also pick a compactly supported time-dependent Hamiltonian function $H$ and an almost complex structure $J$, $\omega$--compatible. When $M$ is not compact but convex at infinity, $J$ will also assumed to be adapted to the boundary (see \cite{McDuffSalamon04} for a detailed description of this condition). From now on, $I$ denotes the interval $[0,1]$.\\

$H$ induces a vector field $X_H\co \bR\times M\ra TM$ given by $\omega(X_H^t,-) = -dH_t$ which gives rise to a family of symplectomorphisms, $\phi_H\co \bR\times M\ra M$, such that $\del_t\phi^t_H = X_H^t\circ\phi_H^t$. We assume that $L'$ and $\phi^1_H(L)$ are transverse. First we define
\begin{align*}
\mathcal{P}(L,L'):=\{ \gamma\in C^\infty(I,M) \, | \; \gamma(0)\in L,\; \gamma(1)\in L'\}
\end{align*}
and we fix $\eta\in \mathcal{P}(L,L')$. Let $\mathcal{P}_\eta (L,L')$ be the connected component of $\mathcal{P}(L,L')$ containing it. The action functional $\mathcal{A}_{H}\co \mathcal{P}_\eta (L,L') \ra \bR$ is given by:
\begin{align*}
\mathcal{A}_H(x)=-\int_{I\times I} \bar{x}^*\omega +\int_I H(t,x(t))\,dt
\end{align*}
where $\bar{x}\co I\times I\ra M$ satisfies the following conditions: $\bar{x}(0,-)=\eta$, $\bar{x}(1,-)=x$ and $\bar{x}(I,0)\subset L$, $\bar{x}(I,1) \subset L'$. Condition (\ref{context}) ensures that the first integral is well-defined.
\begin{rema}\label{remarque:reference action}
In Hamiltonian Floer homology, it is natural to normalize the Hamiltonian. This induces a translation of the action by a constant. In the Lagrangian intersection setting, one can normalize the action functional by assuming that $\mA_H(\eta)=0$. In that case, another choice of the reference path $\eta'$ in the same homotopy class also implies a translation of the action by a constant. Along these lines, we will emphasize our choices of reference as soon as they become crucial.
\end{rema}

The set of the critical points of $\mA_H$, denoted by $\mI(L,L';\eta,H)\subset\mathcal{P}_\eta (L,L')$, consists of orbits of the Hamiltonian vector field. It is in one-to-one correspondence with a subset of $\phi^1_H(L)\cap L'$ (which is finite since we assume transversality and compactness). A Floer trajectory is a smooth map $u\co \bR\times I\ra M$, satisfying
\begin{align*}
\bar{\del}_{J,H}(u):=\del_s u+J_t(u)(\del_t u -X_H^t(u))=0.
\end{align*}
Its energy is defined by the formula
\begin{align*}
E(u) &=\frac{1}{2}\int_I \int_{{\bR}} \left( \|\del_s u\|^2+\|\del_t u - X_H^t(u)\|^2 \right)\,ds dt.
\end{align*}
A Floer trajectory has finite energy if and only if its ends are orbits of the Hamiltonian vector field. For $x$ and $y\in \mI(L,L';\eta,H)$, we form the following moduli space:
\begin{align*}
\widehat{\mathcal{M}}_{x,y}(L,L';H,J) := \left\{ u\in C^\infty(\bR\times I) \, \left|
\begin{array}{l}
\del_s u +J(u)\del_t u + \nabla H_t(u)=0 \\
u(\bR,0)\subset L,\; u(\bR,1)\subset L' \\
u(-\infty,-)=x,\; u(+\infty,-)=y
\end{array}
\right.\right\}.
\end{align*}
All its elements have the same (finite) energy, satisfying
\begin{align}\label{actionenergie}
E(u)= \mA_H(x) - \mA_H(y).
\end{align}
There is a $\bR$--action on $\widehat{\mathcal{M}}_{x,y}(L,L';H,J)$ and we define:
\begin{equation}\label{modhom}
\mM_{x,y}(L,L';H,J) :=\widehat{\mM}_{x,y}(L,L';H,J)/\bR .
\end{equation}

A pair $(H,J)$ is called regular if for all $u\in\widehat{\mathcal{M}}_{x,y}(L,L';H,J)$, the linearization at $u$ of the operator $\bar{\del}_{J,H}$ is surjective. A choice of regular pair $(H,J)$ is generic and for such a pair, $\mM_{x,y}(L,L';H,J)$ is a $(\mu(x)-\mu(y)-1)$--dimensional manifold with $\mu$ the Maslov index. We roughly recall its construction now. For $x$ and $y\in \mI (L,L';\eta,H)$ we choose a smooth map $v\co I\times I\ra M$ such that
\begin{align*}
v(I,0)\subset L,\, v(I,1)\subset L',\, v(0,-)=x\;\; {\rm and}\;\; v(1,-)=y.
\end{align*}
Trivializing $v^*TM$ induces a loop $\gamma_v$ in $\mL (\bR^{2n})$, the set of all Lagrangian submanifolds in $\bR^{2n}$. The degree of $\mu(\gamma_v)$ ($\mu\co \mL (\bR^{2n})\ra S^1$ is the Poincar\'e dual of the Maslov cycle) does not depend on the trivialization or on the choice of $v$ (since $\mu|_{\pi_2(M,L)}=0$). It is denoted $\mu(x,y)$. Since $\mu$ is additive, we fix $z_0$ and define $\mu:=\mu(-,z_0)$.
\begin{rema}\label{remarque:reference maslov}
Choosing another reference $z_0'$ for the Maslov index implies a shift of degrees. Like for the reference of the action functional, this choice will in general be implied and will be emphasized in some key places.
\end{rema}

Similarly to the Morse case, we define
\begin{equation}\label{floercc}
\begin{split}
CF_k(L,L';\eta;H,J) & := \langle x\in \mI(L,L';\eta,H) |\; \mu(x)=k\rangle_{\bZ_2} , \\
\del(x) & :=\sum_{y|\, \mu(x,y)=1} \#_2 \mathcal{M}_{x,y}(L,L';H,J) \cdot y.
\end{split}
\end{equation}

When $\mu(x,y)\neq 1$, $\mM_{x,y}(L,L';H,J)$ is not compact in general but admits a compactification (Gromov's compactness theorem for the compact or convex at infinity setting). As condition (\ref{context}) prohibits bubbling we get (the case $\mu(x,y)>1$ has been treated by Barraud and Cornea \cite{BarraudCornea07})
\begin{align*}
\del\overline{\mM}_{x,y}(L,L';H,J)= \bigcup_{z\in \mI(L,L';\eta,H)} \overline{\mM}_{x,z}(L,L';H,J)\times \overline{\mM}_{z,y}(L,L';H,J).
\end{align*}
In particular, this implies that $\del^2 =0$ and thus (\ref{floercc}) defines a chain complex. Its homology is the Lagrangian Floer homology of $M$ relatively to $L$ and $L'$:
\begin{align*}
HF_*(L,L'):=H(CF_*(L,L';\eta;H,J),\del).
\end{align*}
It is well-known that this homology does not depend on the choice of the regular pair $(H,J)$ (see \S\ref{subsec:comparison morphisms} and \S\ref{subsec:pss} for two different proofs of this fact).


\subsection{Comparison morphisms}

We present here three morphisms which all induce an isomorphism in homology.
\begin{itemize}
\item $\psi^{01}\co CF_*(L,L';\eta;H_0,J_0)\ra CF_*(L,L';\eta;H_1,J_1)$ \hfill(comparison morphism)
\item $b_H\co CF_*(L,L';\eta;H,J)\ra CF_*(L,L'';\widetilde{\eta};0,\widetilde{J})$ \hfill(naturality isomorphism)
\item $\phi_f^{H}\co CM_*(L;f,g)\ra CF_*(L,L;\eta;H,J)$ \hfill(Lagrangian PSS morphism)
\end{itemize}

\subsubsection{Classical comparison morphism of Lagrangian Floer homology}\label{subsec:comparison morphisms}

This morphism compares the Floer complexes, for two different regular pairs $(H_0,J_0)$ and $(H_1,J_1)$. It induces an isomorphism in homology.

Let $H^{01}\co\bR\times ([0,1]\times M)\ra \bR$ and $J^{01}\co\bR\times M\ra {\rm End}(TM)$ be respectively $1$--parameter families of Hamiltonian functions and almost complex structures $\omega$--compatible. We assume that there exists $R>0$ such that for all $s\leq -R$, $(H_s^{01},J_s^{01})=(H_0,J_0)$ and all $s\geq R$, $(H_s^{01},J_s^{01})=(H_1,J_1)$ ($H_s^{01}$ and $J_s^{01}$ denote $H^{01}(s,-,-)$ and $J^{01}(s,-)$). Moreover we assume that there is a compact set, containing the support of all these Hamiltonian functions.\\

$H^{01}$ induces a family of vector fields $X^{01}$ and a family of symplectomorphisms $\phi^{01}$ via $\omega(X^{01},-) = -dH^{01}$ and $\del_t\phi^{01} = X^{01}(\phi^{01})$. The equation satisfied by the Floer trajectories $u\co\bR\times [0,1]\ra M$ becomes
\begin{align*}
\del_s u(s,t) +J^{01}_s(u(s,t))\del_t u(s,t) + \nabla^s_x H_s^{01}(t,u(s,t)) =0.
\end{align*}
$\nabla^s_x H^{01}$ is the gradient of $H^{01}$ induced by the metric associated to $J_s^{01}$. The boundary conditions are $u(\bR,0)\subset L$ and $u(\bR,1)\subset L'$. We also require their energy to be finite:
\begin{align*}
E(u)= \frac{1}{2}\int_{{\bR}\times I} \| \del_s u\|^2 + \| \del_t u - X^{01}(u)\|^2 dsdt < \infty
\end{align*}
which is anew equivalent to the existence of orbits $x_0\in\mI(L,L';\eta,H_0)$ and $y_1\in\mI(L,L';\eta,H_1)$ such that $u(-\infty,-)=x_0$ and $u(+\infty,-)=y_1$.
For generic choice $(H_s^{01},J_s^{01})$ is regular and for such a regular pair $\mM_{x_0,y_1}(L,L';H^{01},J^{01})$ is a smooth $\mu(x_0,y_1)$--dimensional manifold. After compactification, the boundary of these moduli spaces is the disjoint union:
\begin{align*}
&\bigcup_{y_0\in \mI(L,L';\eta,H_0)} {\overline{\mM}}_{x_0,y_0}(L,L';H_{0},J_{0}) \times {\overline{\mM}}_{y_0,y_1}(L,L';H^{01},J^{01})\\
&\bigcup_{x_1\in\mI(L,L';\eta,H_1)} {\overline{\mM}}_{x_0,x_1}(L,L';H^{01},J^{01}) \times {\overline{\mM}}_{x_1,y_1}(L,L';H_{1},J_{1}).
\end{align*}
Therefore, the formula
\begin{align*}
\psi^{01}(x_0):= \sum_{y_1\in\mI (L,L';\eta,H_1) |\, \mu(x_0,y_1)=1}  \#\mM_{x_0,y_1}(L,L';H^{01},J^{01})\cdot y_1
\end{align*}
is easily seen to define a morphism of chain complexes inducing an isomorphism $\psi^{01}\co HF_*(L,L';H_0,J_0)\ra HF_*(L,L';H_1,J_1)$.

\subsubsection{Naturality}

This morphism appears for example in Barraud and Cornea~\cite{BarraudCornea07}. It is an action preserving identification of the complexes $CF_*(L,L';\eta;H,J)$ and $CF_*(L,L'';\widetilde{\eta};0,\widetilde{J})$, where $L'':= (\phi^1_H)^{-1}(L')$, $\widetilde{J}:= \phi^* J$ and $\eta$, $\widetilde{\eta}$ satisfying $\eta(t)=\phi^t_H(\widetilde{\eta}(t))$.\\

First we define the map
\begin{align*}
b_H\co \left\{
\begin{array}{rcl}
\mP(L,L'') &\ra &\mP(L,L')\\
x &\mapsto &[t\mapsto \phi^t_H(x(t))]
\end{array}
\right.
\end{align*}
It induces a map $b_H\co \mP_{\widetilde{\eta}}(L,L'')\ra\mP_\eta (L,L')$. Since $L \cap L''$ is transverse, $b_H$ identifies $\mI(L,L'';\widetilde{\eta},0)$ and $\mI(L,L';\eta,H)$.\\

For $u\in\mM_{x,y}(L,L'';0,\widetilde{J})$, let $\widetilde{u}$ be given by $\widetilde{u}(s,t):= \phi^t_H(u(s,t))$. Comparing the linearization of the two operators $\bar{\del}_{J,H}$ and $\bar{\del}_{\widetilde{J},0}$ ensures that $(H,J)$ is a regular pair if and only if $(0,\widetilde{J})$ is. One can show that $\bar{\del}_{J,H}(\widetilde{u}) = (\phi_H)_*(\bar{\del}_{\widetilde{J},0}(u))$ for any $u\in\mM_{x,y}(L,L';H,J)$. Moreover, boundary conditions correspond. Since we have
\begin{align*}
\forall x\in \mI(L,L'';\eta',0), \;\; \mA_{H}(b_H(x))=\mA_{0}(x)
\end{align*}
we get, via (\ref{actionenergie}), $E(\widetilde{u})=E(u)$. Notice that the two action functionals correspond for our choices of reference $\eta$ and $\widetilde{\eta}$ if we assume that they have been normalized (see Remark \ref{remarque:reference action}). Therefore $b_H$ induces a bijection
\begin{align*}
b_H\co\mM_{x,y}(L,L'';0,\widetilde{J})\ra \mM_{b_H(x),b_H(y)}(L,L';H,J)
\end{align*}
which gives rise to an action preserving identification of chain complexes:
\begin{align*}
b_H\co  CF_*(L,(\phi^1_H)^{-1}(L');\widetilde{\eta};0,\phi^*J)\ra CF_*(L,L';\eta;H,J)
\end{align*}

\subsubsection{Lagrangian PSS morphism}\label{subsec:pss}

The Hamiltonian Piunikhin--Salamon--Schwarz (PSS) morphism \cite{PiunikhinSalamonSchwarz96} has been adapted to the Lagrangian setting, for the particular case of cotangent bundles by Kati\'c and Milinkovi\'c \cite{KaticMilinkovic05} and in more generality by Barraud and Cornea~\cite{BarraudCornea05}, and Albers~\cite{Albers05}. This morphism compares the Morse and Floer complexes and induces an isomorphism in homology.

\begin{figure}[h]\label{defPSS}
\centering
\psfrag{B1}[][][1]{\footnotesize(\ref{pssbord1})}
\psfrag{B2}[][][1]{\footnotesize(\ref{pssbord2})}
\psfrag{g(0)}[][][1]{\footnotesize$\gamma(0)$}
\psfrag{g}[][][1]{\footnotesize$\gamma$}
\psfrag{p}[][][1]{\footnotesize$p$}
\psfrag{p'}[][][1]{\footnotesize$p'$}
\psfrag{u}[][][1]{\footnotesize$u$}
\psfrag{x}[][][1]{\footnotesize$x$}
\psfrag{x'}[][][1]{\footnotesize$x'$}
\psfrag{L}[][][1]{\footnotesize$L$}
\includegraphics{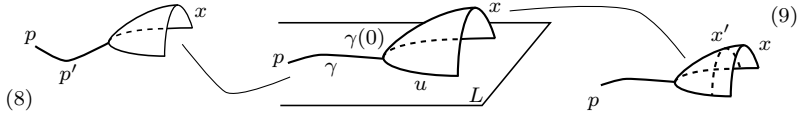}
\caption{Definition of the Lagrangian PSS morphism}
\label{PSSlagr}
\end{figure}

To $p\in {\rm Crit}(f)$ and $x\in \mI(L,L;\eta,H)$, we associate
\begin{align*}
\mM^f_{p}(g) &:= \left\{ \gamma:\bR\ra L\left|\, - \nabla f(\gamma(t)) = \frac{d\gamma(t)}{dt}, \, \gamma(-\infty)=p \right. \right\} =\mW^u_p(f,g),\\
\mM^H_{x}(J) &:= \left\{ u\co\bR\times I\ra M \left|
\begin{array}{l}
\del_s u+J\del_t u +\beta(s)\nabla H(u)=0\\
u(+\infty,t)=x(t), u(\bR,\{ 0,1 \})\subset L
\end{array}
 \right. \right\}
\end{align*}
where $\beta(s)$ is a smooth, increasing function, whose value is $0$ for $s\leq 1/2$ and $1$ for $s\geq 1$. Up to generic choices, $\mM^H_{x}(J)$ is a smooth $\mu(x)$--dimensional manifold and
\begin{align*}
\mM^{f,H}_{p,x}:= \left\{ \left. (\gamma,u)\in \mM^f_{p}(g)\times \mM^H_{x}(J) \right| \, u(-\infty,-)=\gamma(0) \right\}
\end{align*}
a ($i_f(p)-\mu(x)$)--dimensional manifold (notice that this requires a particular choice of the reference of the Maslov index). Its $0$--dimensional component is compact and its $1$--dimensional component admits a compactification such that (see Figure \ref{PSSlagr})
\begin{align}
\del \overline{\mM}^{f,H}_{p,x} = &\,\;\;\bigcup_{p'\in {\rm Crit}(f)} \,\;\;\mM_{p,p'}(f,g) \times \mM^{f,H}_{p',x}\label{pssbord1}\\
&\bigcup_{x'\in\mI(L,L;\eta,H)} \mM^{f,H}_{p,x'} \times \mM_{x',x}(L,L;H,J).\label{pssbord2}
\end{align}

\begin{figure}[h]
\centering
\psfrag{g}[][][1]{\footnotesize$\gamma$}
\psfrag{g+}[][][1]{\footnotesize$\gamma^+$}
\psfrag{g-}[][][1]{\footnotesize$\gamma^-$}
\psfrag{g'}[][][1]{\footnotesize$\gamma_1$}
\psfrag{g0}[][][1]{\footnotesize$\gamma_0$}
\psfrag{x}[][][1]{\footnotesize$x$}
\psfrag{*}[][][1]{\footnotesize$\ast$}
\psfrag{y'}[][][1]{\footnotesize$y_1$}
\psfrag{y}[][][1]{\footnotesize$y_0$}
\psfrag{q}[][][1]{\footnotesize$q$}
\psfrag{p}[][][1]{\footnotesize$p_0$}
\psfrag{p'}[][][1]{\footnotesize$p_1$}
\psfrag{R}[][][1]{\footnotesize$R$}
\psfrag{R>0}[][][1]{\footnotesize$R\ra 0$}
\psfrag{R>infty}[][][1]{\footnotesize$R\ra \infty$}
\psfrag{u0}[][][1]{\footnotesize$u_0$}
\psfrag{u-}[][][1]{\footnotesize$u-$}
\psfrag{u+}[][][1]{\footnotesize$u+$}
\psfrag{u'}[][][1]{\footnotesize$u_1$}
\psfrag{u}[][][1]{\footnotesize$u$}
\psfrag{unicité}[][][1]{\footnotesize unicity}
\psfrag{gluing}[][][1]{\footnotesize gluing}
\includegraphics{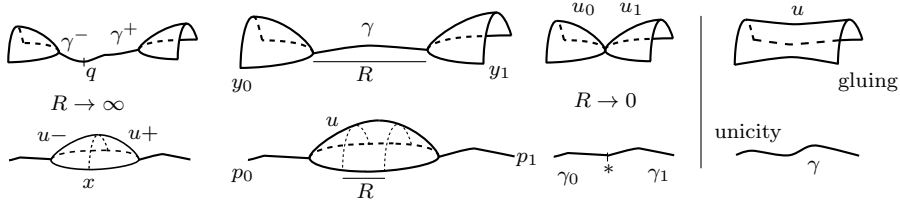}
\caption{Commutativity of the PSS and comparison morphisms}
\label{commPSScomp}
\end{figure}

The morphism $\phi_f^{H}\co CM_*(L;f,g)\ra CF_*(L,L;H,J)$, defined on generators by
\begin{align*}
\phi_f^{H}(p):=\sum_{x|\, \mu(x)=i_f(p)} \#_2 \mM^{f,H}_{p,x}\cdot x
\end{align*}
induces a morphism $\phi_f^{H}\co HM_*(L;f,g)\ra HF_*(L,L;H,J)$. Similarly, the moduli spaces
\begin{align*}
\mM_f^{q}(g) &:= \left\{ \gamma\co\bR\ra L\left|\, - \nabla f(\gamma(t)) = \frac{d\gamma(t)}{dt}, \, \gamma(+\infty)=q \right. \right\} = \mW^s_q(f,g),\\
\mM_H^{y}(J) &:= \left\{ u\co\bR\times I\ra M \left|
\begin{array}{l}
\del_s u+J\del_t u +\beta(-s)\nabla H(u)=0\\
u(-\infty,t)=y(t), u(\bR,\{ 0,1 \})\subset L
\end{array}
\right. \right\} \;\;{\rm and}\\
\mM_{H,f}^{y,q} &:= \left\{ \left. (u,\gamma)\in \mM_H^{y}(J)\times \mM_f^{q}(g) \right| \, \gamma(0)=u(+\infty,-) \right\}
\end{align*}
allow us to define a morphism
\begin{align*}
\psi_{H}^f(y):=\sum_{q|\, i_f(q)=\mu(y)} \#_2 \mM_{H,f}^{y,q}\cdot q
\end{align*}
which induces $\psi^f_{H}\co HF_*(L,L;H,J)\ra HM_*(L;f,g)$. These two morphisms commute with the classical comparison morphisms, as stated below.
\begin{lemm}\label{lemm:PSScomparaison}
The following diagrams commute:
\begin{align}
 \begin{split}
\xymatrix{\relax
    HF_*(L,L;H_0,J_0) \ar[rd]_{\psi_{H_0}^f}\ar[rr]^{\psi^{01}} && HF_*(L,L;H_1,J_1)\\
    & HM_*(L;f,g) \ar[ru]_{\phi^{H_1}_f}}\label{psscomparaison}
\end{split}\\
 \begin{split}
\xymatrix{\relax
    HM_*(L;f_0,g_0) \ar[rd]_{\phi^H_{f_0}} \ar[rr]^{\psi^{01}_{\rm Morse}} && HM_*(L;f_1,g_1)\\
    & HF_*(L,L;H,J) \ar[ru]_{\psi_H^{f_1}}}\label{psscomparaisonmorse}
\end{split}
\end{align}
\end{lemm}
To prove it, we study cobordisms (see Figure \ref{commPSScomp}) similar to the ones giving that $\phi_f^{H}\circ\psi_H^f$ and $\psi^f_{H}\circ\phi^H_f$ induce the identity in homology (which is an immediate corollary of the lemma, for $(H_0,J_0) = (H_1,J_1)$ and $(f_0,g_0)=(f_1,g_1)$).

\begin{proof} {\it Diagram (\ref{psscomparaison})} -- Let $\{ J_s^{01} \}_{s\in {\Bbb R}}$ be a $1$--parameter family of almost complex structures, such that $J^{01}_s=J_0$, if $s\leq -1$, $J^{01}_1=J$ if $|s|<1/2$ and $J^{01}_1=J_1$ for $s\geq 1$. To $y_i\in \mI(L,L;\eta,H_i)$ ($i=0,1$), we define moduli spaces as follows:
\begin{multline*}
\mM^{\phi\circ\psi}_R(y_0,y_1):= \\
\left\{ (u_0,\gamma,u_1) \in \mM^{y_0}_{H_0}(\bar{J}_0)\times C^{\infty}([0,R],M)\times \mM_{y_1}^{H_1}(\bar{J}_1) \left|
\begin{array}{l}
\nabla f(\gamma(t)) = \frac{d\gamma(t)}{dt}\\
u_0(+\infty)=\gamma(0)\\
u_1(-\infty)=\gamma(R)
\end{array}
\right. \right\}
\end{multline*}
where $\bar{J}_0$ and $\bar{J}_1$ are $1$--parameter families of almost complex structures, such that
\begin{align*}
\bar{J}_0=J_0 \mbox{ if } s\leq -1 &\mbox{  and } \bar{J}_0=J \mbox{ if } s\leq -1/2,\\
\bar{J}_1=J_1 \mbox{ if } s\geq 1 &\mbox{  and } \bar{J}_1=J \mbox{ if } s\geq 1/2.
\end{align*}
We define $\mM^{\phi\circ\psi}(y_0,y_1):= \bigcup_{R>0} \mM^{\phi\circ\psi}_R(y_0,y_1)$. Up to generic choices, it is a smooth ($\mu(y_0)-\mu(y_1)+1$)--dimensional manifold. Its $0$--dimensional component is compact. Hence one can define a chain homotopy $\varphi\co CF_*(L,L;H_0,J_0) \ra CF_*(L,L;H_1,J_1)$, by the formula
\begin{align*}
\varphi(y_0):=\sum_{y_1|\, \mu(y_1)=\mu(y_0)+1} \#_2 \mM^{\phi\circ\psi}(y_0,y_1)\cdot y_1.
\end{align*}
Its $1$--dimensional component admits a compactification such that
\begin{align}
\del\overline{\mM}^{\phi\circ\psi}(y_0,y_1) = &\bigcup_{x_0\in \mI(L,L;\eta,H_0)} \mM_{y_0,x_0}(L,L';H_0,J_0)\times \mM^{\phi\circ\psi}(x_0,y_1)\label{psscomm1}\\
&\;\;\;\;\;\;\;\,\;\bigcup \;\;\;\;\;\;\;\mM^{\phi\circ\psi}_{0}(y_0,y_1)\label{psscomm2}\\
&\;\;\;\;\,\bigcup_{q\in {\rm Crit}(f)} \;\;\;\mM_{H_0,f}^{y_0,q}\times \mM^{f,H_1}_{q,y_1}\label{psscomm3}\\
&\,\bigcup_{x_1\in \mI(L,L;\eta,H_1)} \mM^{\phi\circ\psi}(y_0,x_1)\times \mM_{x_1,y_1}(L,L';H_1,J_1)\label{psscomm4}
\end{align}
(\ref{psscomm1}) and (\ref{psscomm4}) arise when $R$ converges to a real number, (\ref{psscomm2}), (\ref{psscomm3}) when it converges to $0$ and to infinity. Boundary (\ref{psscomm2}), $\mM^{\phi\circ\psi}_0(y_0,y_1)$, consists of pairs of discs which are pseudo-holomorphic near their common point. They can be glued (see~\cite{FukayaOhOhtaOno} and~\cite{BiranCornea}) to give the Floer trajectories defining the comparison morphism and then $\del\varphi + \varphi\del = \phi^{H_1}_f\circ \psi_{H_0}^f + \psi^{01}$. Diagram (\ref{psscomparaison}) commutes.\\

{\it Diagram (\ref{psscomparaisonmorse})} -- Given two pairs $(f_0,g_0)$, $(f_1,g_1)$, let $(f^{01},g^{01})$ be a regular homotopy such that $(f_s^{01},g_s^{01})=(f_0,g_0)$ when $s<-1$, $(f_s^{01},g_s^{01}) =(f_1,g_1)$ when $s>1$. We define the sets
\begin{align*}
\mM_R^{\psi\circ\phi}(p_0,p_1):= {\bigg \{} (\gamma_0,U_R,\gamma_1) \in \mW^u_{p_0}(f^{01},g^{01})\times C^{\infty}(\bR\times I,M)\times \mW^s_{p_1}(f^{01},g^{01}) \\
\left|
\begin{array}{l}
\del_s U_R + J\del_t U_R + \alpha_R(s) \nabla H(U_R)=0,\\
U_R(-\infty,-)=\gamma_0(0), \; U_R(+\infty,-)=\gamma_1(0)
\end{array}
\right. {\bigg \}}
\end{align*}
where $\alpha_R$ is a smooth cut-off function, whose value is $1$ for $|s|<R$ and $0$ for $|s|>R+1$. We also assume that it is $C^1$--bounded to get an estimate on discs energy -- see~\cite{Albers05}). Let define $\mM^{\psi\circ\phi}(p_0,p_1):= \bigcup_{R>0} \mM_R^{\psi\circ\phi}(p_0,p_1)$ which is a smooth ($i_{f^{01}}(p_0)-i_{f^{01}}(p_1)+1$)--dimensional manifold for generic choices. The boundary of its (compactified) $1$--dimensional component is
\begin{align}
\del\overline{\mM}^{\psi\circ\phi}(p_0,p_1) = &\;\bigcup_{{q_0}\in {\rm Crit}(f^{01})} \;\mM_{p_0,q_0}(f^{01},g^{01})\times \mM_R^{\psi\circ\phi}(q_0,p_1)\label{psscomm1'}\\
&\;\;\;\;\;\;\;\bigcup \;\;\;\;\;\;\;\mM^{\psi\circ\phi}_{0}(p_0,p_1)\label{psscomm2'}\\
&\bigcup_{x\in \mI(L,L;\eta,H)} \mM^{f_0,H}_{p_0,x}\times \mM_{H,f_1}^{x,p_1}\label{psscomm3'}\\
&\;\bigcup_{{q_1}\in {\rm Crit}(f^{01})} \;\mM_R^{\psi\circ\phi}(p_0,q_1)\times \mM_{q_1,q_0}(f^{01},g^{01})\label{psscomm4'}
\end{align}
The parameter $R$ converges to a real number in (\ref{psscomm1'}), (\ref{psscomm4'}), to infinity in (\ref{psscomm3'}) and to $0$ in (\ref{psscomm2'}). When $R$ is infinite, the products lead to the composition $\phi^{H}_{f_0}\circ\psi_{H}^{f_1}$, whereas when $R=0$, one gets a pseudo-holomorphic disc with boundary in $L$, with null symplectic area (since $\omega|_{\pi_2(M,L)}=0$). Hence it has to be constant. Therefore, when $R$ converges to $0$, we get two critical points $p_0$ and $p_1$, with a Morse flow line (of $f^{01}$) passing through them (flow lines unicity), with the same indices. Thus we recover the moduli spaces defining $\psi^{01}_{\rm Morse}$, which proves that diagram (\ref{psscomparaisonmorse}) commutes.
\end{proof}


\section{Proof of Theorem \ref{theo:independ}}\label{sec:demo}

We consider the action filtration of the Floer complex. For a real number $\nu$ and an integer $k$, one can define the $\bZ_2$--vector space:
\begin{align*}
CF_k^\nu (L,L';\eta;H,J) := \langle x\in \mI(L,L';\eta,H) \, | \;\mu(x)=k, \; \mA_{H}(x)\leq \nu \rangle_{\bZ_2}.
\end{align*}
As the action functional decreases along Floer trajectories, Floer's differential preserves this filtration and $(CF_k^\nu (L,L';\eta;H,J),\del)$ is a subcomplex. We denote its inclusion by $i^\nu\co CF_*^\nu(L,L';\eta;H,J)\ra CF_*(L,L';\eta;H,J)$ and by $i^\nu_*:= H_*(i^\nu)$. We define for $x\in HF_*(L,L';\eta;H,J)$, $x\neq 0$
\begin{align*}
\widetilde{\sigma}_{L,L'}(x;H,J,\eta) := \inf \{ \nu\in\bR \, |\; x\in {\rm Im}(i_*^\nu) \}
\end{align*}
and we put for $\alpha\in HM_*(L;f,g)$, $\alpha\neq 0$, (see Definition \ref{defi:si})
\begin{align*}
\sigma_{L}(\alpha;H,J,\eta;f,g) &:= \widetilde{\sigma}_{L,L}(\phi^H_f(\alpha);H,J,\eta)\;\;{\rm and}\\
c(\alpha;H,J;f,g) &:= \sigma_{L}(\alpha;H,J,\eta;f,g) - \sigma_{L}(1;H,J,\eta;f,g)
\end{align*}
with $1$ the generator of $HM_0(L;f,g)$. By definition, homological Lagrangian spectral numbers are critical values of the action functional. The commutativity of diagram (\ref{psscomparaisonmorse}) gives their independence on the Morse--Smale pair $(f,g)$. Hence, we will denote them by $\sigma_L(\alpha;H,J,\eta)$ et $c_L(\alpha;H,J)$ (it is clear in view of Remark \ref{remarque:reference action} that the absolute version $c_L(-;H,J)$ does not depend on the choice of $\eta$).\\

The aim of this section is to prove Theorem \ref{theo:independ}. The most delicate part follows from the proposition below, which is also of interest in itself.

Indeed this proposition states, in particular, that, given a Morse--Smale pair $(f,g)$ and a regular pair $(H,J)$, the composition of the PSS morphism and (the inverse of) the naturality morphism
\begin{align*}
b_H^{-1}\circ \phi_f^H \co HM_*(L;f,g)\lra HF_*(L,(\phi_H^1)^{-1}(L);0,\phi_H^*J)
\end{align*}
only depends on $(\phi_H^1)^{-1}(L)$ and $\phi_H^*J$. As the naturality morphism preserves the action, $b_H^{-1}\circ \phi_f^H$ is, from this point of view, a more natural candidate than the PSS morphism in the study of spectral numbers.

\begin{prop}\label{prop:commut}
For every regular pairs $(H,J)$ and $(H',J')$, such that
\begin{align*}
\phi_H^*J=\phi_{H'}^*J'=:\widetilde{J} \;\mbox{ and }\; (\phi_H^1)^{-1}(L)=(\phi^1_{H'})^{-1}(L)=:L_0
\end{align*}
the following diagram commutes:
\begin{align}\label{diag}
 \begin{split}
\xymatrix{\relax
    HM_*(L;f,g) \ar[d]_{\phi_f^{H'}}\ar[r]^{\hspace{-.3cm}\phi_f^{H}} & HF_*(L,L;H,J) \ar[d]^{b_H^{-1}}\\
    HF_*(L,L;H',J') \ar[r]_{b_{H'}^{-1}} & HF_*(L,L_0;0,\widetilde{J})
}
 \end{split}
\end{align}
\end{prop}

\hop{\it Organization of the next four sections} -- In the next two, we discuss additional algebraic structures required to prove Proposition \ref{prop:commut}. In Section \ref{sec:structure}, we show that Floer homology can be viewed as a module over Morse homology (endowed with the intersection product). In Section \ref{sec:independance}, we show that the naturality and the Lagrangian PSS morphisms agree with this structure which will be seen to immediately prove Proposition \ref{prop:commut}. In Section \ref{sec:continuity}, we show the independence on the almost complex structure and the continuity property. Then we collect all the results such as to prove Theorem \ref{theo:independ}. In Section \ref{sec:classical}, we prove an additional property, namely that these invariants are a generalization of the Hamiltonian ones, by showing that equality (\ref{prop:generalisation}) holds.


\subsection{Lagrangian Floer homology as a module over Morse homology}\label{sec:structure}

In this section, we recall that $HM_*(L)$ is a unitary ring and we equip $HF_*(L,L')$ with a module structure over it.

\subsubsection{Morse homology ring}\label{section:homology de Morse anneau}

We recall the Morse theoretic version of the intersection product. It is defined at the chain level:
\begin{align*}
\xymatrix{
CM_k(L;f_1,g) \otimes CM_l(L;f_2,g)\ar[r]^{\hspace{.85cm}\cdot} & CM_{k+l-n}(L;f_3,g) }
\end{align*}
where $f_1$, $f_2$ and $f_3$ are Morse functions and $g$ a metric such that for all $i$ and $j\in \{ 1,2,3\}$, $\mW_{p_i}^u(f_i,g)$ and $\mW_{p_j}^s(f_j,g)$ intersect transversally for any $p_i$, $p_j$ respectively critical points of $f_i$ and $f_j$. This implies in particular that the three pairs $(f_i,g)$ ($i=1,3$) are Morse--Smale. These choices are generic.\\

To $p\in {\rm Crit}(f_1)$, $q\in {\rm Crit}(f_2)$ and $r\in {\rm Crit}(f_3)$, we associate the moduli space:
\begin{multline*}
\mM_{(p,q);r}(f_1,f_2,f_3;g):= \\
 \left\{ (\gamma_1,\gamma_2,\gamma_3) \in \mW^u_p(f_1,g)\times \mW^u_q(f_2,g)\times \mW^s_r(f_3,g) |\, \gamma_1(0)=\gamma_2(0)=\gamma_3(0)
\right\}.
\end{multline*}
Under our assumptions, these spaces are manifolds of dimension:
\begin{align*}
{\rm d}(p,q;r)= i_{f_1}(p)+i_{f_2}(q)-i_{f_3}(r)-n.
\end{align*}

Their $0$--dimensional component is compact, then one can define the product on the generators of the chain complex (and extend it by bilinearity):
\begin{align*}
p\cdot q := \sum_{r|\,{\rm d}(p,q;r)=0} \#_2 \mM_{(p,q);r}(f_1,f_2,f_3;g)\cdot r.
\end{align*}
One can show that this formula induces a product in homology which is independent on the choices of metric and Morse functions. We recall that the unity of $(HM_*(L),\cdot)$ is $[L]$, the fundamental class of $L$.

\subsubsection{Floer homology as a module over Morse homology}

We define, analogously, an external product on $HF_*(L,L')$:
\begin{align*}
\xymatrix{
HM_*(L) \otimes HF_*(L,L')\ar[r]^{\hspace{.6cm}\ast} & HF_*(L,L') }.
\end{align*}
As before, it is defined at the chain level. In order to shorten notation, $L$ and $L'$ do not appear in what follows. For example, the moduli spaces $\mM_{x,y}(L,L';H,J)$ defined by (\ref{modhom}) are denoted here and in Section \ref{sec:independance} by $\mM_{x,y}(H,J)$.

\begin{figure}[h]
\centering
\psfrag{B1}[][][1]{\footnotesize(\ref{pdtfloerbord1})}
\psfrag{B2}[][][1]{\footnotesize(\ref{pdtfloerbord2})}
\psfrag{B3}[][][1]{\footnotesize(\ref{pdtfloerbord3})}
\psfrag{p}[][][1]{\footnotesize$p$}
\psfrag{x}[][][1]{\footnotesize$x$}
\psfrag{y}[][][1]{\footnotesize$y$}
\psfrag{p'}[][][1]{\footnotesize$p'$}
\psfrag{x'}[][][1]{\footnotesize$x'$}
\psfrag{y'}[][][1]{\footnotesize$y'$}
\psfrag{g}[][][1]{\footnotesize$\gamma$}
\psfrag{u}[][][1]{\footnotesize$u$}
\includegraphics{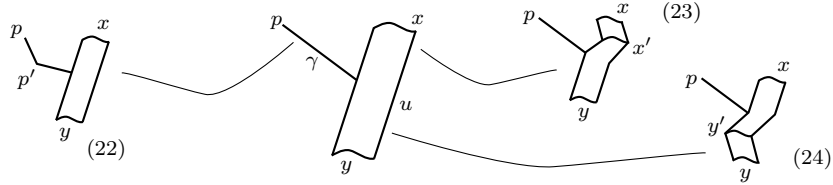}
\caption{Definition of the module structure}
\label{modulestr}
\end{figure}

We choose a Morse--Smale pair $(f,g)$ and a regular pair $(H,J)$. We associate, to $p\in {\rm Crit}(f)$ and $x$, $y\in\mI(L,L';\eta,H)$, the following moduli space:
\begin{align}\label{modpdt}
\mM_{(p,x);y}(f,g;H,J):= \left\{ (\gamma,u)\in \mW^u_p(f,g)\times \mM_{x,y}(H,J) |\, \gamma(0)=u(0,0) \right\}.
\end{align}
These spaces are manifolds of dimension ${\rm d}(p,x;y)= \mu(x)-\mu(y)+i_f(p)-n$ as preimage of the diagonal $\Delta\subset L\times L$ by the map
\begin{align*}
({\rm ev},{\rm ev})\co \left\{
\begin{array}{rcl}
\mW^u_p(f,g)\times \mM_{x,y}(H,J) &\ra &L\times L\\
(\gamma,u) &\mapsto &(\gamma(0),u(0,0))
\end{array}
\right.
\end{align*}
which generically is transverse to $\Delta$. The $0$--dimensional component of these spaces is compact, thus the formula
\begin{align*}
p\ast x:= \sum_{y|\, {\rm d}(p,x;y)=0} \#_2 \mM_{(p,x);y}(f,g;H,J) \cdot y
\end{align*}
defines a product $\ast\co CM_k(L;f,g)\times CF_l(L,L';H,J) \ra CF_{k+l-n}(L,L';H,J)$, when extended by bilinearity. Moreover, the boundary of the compactification of their $1$--dimensional component is the disjoint union of (see Figure \ref{modulestr})
\begin{align}
&\;\;\;\bigcup_{p'\in{\rm Crit}(f)} \;\;\;\mM_{p,p'}(f,g) \times \mM_{(p',x);y}(f,g;H,J)\label{pdtfloerbord1}\\
& \bigcup_{x'\in\mI(L,L';\eta,H)} \mM_{x,x'}(H,J) \times \mM_{(p,x');y}(f,g;H,J)\label{pdtfloerbord2}\\
& \bigcup_{y'\in\mI(L,L';\eta,H)} \mM_{(p,x);y'}(f,g;H,J) \times \mM_{y',y}(H,J)\label{pdtfloerbord3}
\end{align}

Hence, $\ast$ commutes with the respective boundary maps and induces a product in homology.

\begin{lemm}\label{floermodule}
$(HF_*(L,L';H,J),\ast)$ is a $HM_*(L;f,g)$--module.
\end{lemm}

\begin{proof}
The main difficulty is to prove that, for any $\alpha$, $\beta$ in $HM_*(L;f,g)$ and any $a$ in $HF_*(L,L';H,J)$, the following equality holds:
\begin{align*}
(\alpha\cdot \beta)\ast a = \alpha\ast (\beta\ast a).
\end{align*}
We fix two Morse functions $f_1$ and $f_2$ and a metric $g$ satisfying the extended Morse--Smale condition defined at the beginning of \S\ref{section:homology de Morse anneau}. We also choose a regular pair $(H,J)$. We have to prove that $\varphi^\pm$, respectively defined by $\varphi^+ (p,q,x):=(p\cdot q)\ast x$ and $\varphi^-(p,q,x):= p\ast (q\ast x)$ (for $p\in {\rm Crit}(f_1)$, $q\in {\rm Crit}(f_2)$ and $x\in \mI(L,L';\eta,H)$), induce the same morphism in homology. We define respectively $\mM^{+}_{(p,q,x);y}(f_1,f_2,g;H,J)$ and $\mM^{-}_{(p,q,x);y}(f_1,f_2,g;H,J)$ by
\begin{align*}
&\,\;\;\bigcup_{r\in{\rm Crit}(f_1)} \,\;\;\mM_{(p,q);r}(f_1,f_2,f_1;g) \times \mM_{(r,x);y}(f_1,g;H,J) \;\; {\rm and}\\
&\bigcup_{z\in\mI(L,L';\eta,H)} \mM_{(q,x);z}(f_2,g;H,J) \times \mM_{(p,z);y}(f_1,g;H,J).
\end{align*}
We recall that $\mM_{(p,q);r}(f_1,f_2,f_1;g)$ denotes the moduli spaces used in order to define the module structure on $HM_*(L)$ (\S\ref{section:homology de Morse anneau}). These spaces are such that
\begin{align*}
\varphi^\pm (p,q,x) =\sum_y \#_2 \mM^\pm_{(p,q,x);y}(f_1,f_2,g;H,J)\cdot y.
\end{align*}
In the following two steps, we define moduli spaces $\mM^0_{(p,q,x);y}(f_1,f_2,g;H,J)$ and a morphism
\begin{align*}
\varphi^{0}(p,q,x):= \sum_y \#_2 \mM^{0}_{(p,q,x);y}(f_1,f_2,g;H,J)\cdot y
\end{align*}
and we conclude the proof by showing that
\begin{align}
H_*(\varphi^+) &= H_*(\varphi^0) = H_*(\varphi^-).\label{module}
\end{align}
Notice that, when $q= m$, the maximum of $f_2$, $i_{f_2}(m)=n$ and thus the moduli spaces defining $m\ast x$ and $p\cdot m$ have such a dimension that, even on the chain level, we have $m\ast x =x$ and $p\cdot m=p$. Hence, it holds
\begin{align}\label{rema:max}
\varphi^+ (p,m,x)=(p\cdot m)\ast x= p\ast x = p\ast (m\ast x)=\varphi^- (p,m,x).
\end{align}

\hop {\it (step 1)} First, for $R>0$, we consider the moduli spaces
\begin{multline*}
\mM^R_{(p,q,x);y}(f_1,f_2,g;H,J):= \\
\left\{ (\gamma^R_1,\gamma_2,u)\in \mW^u_p(f_1,g)\times \mW^u_q(f_2,g)\times \mM_{x,y}(H,J) \left|
\begin{array}{l}
\gamma^R_1(0)=\gamma_2(0) \\
\gamma^R_1(R)=u(0,0)
\end{array}
\right. \right\}
\end{multline*}
where $p\in {\rm Crit}(f_1)$, $q\in {\rm Crit}(f_2)$, $x$ and $y\in\mI(L,L';\eta,H)$. We define
\begin{align*}
\mM_{(p,q,x);y}(f_1,f_2,g;H,J):= \bigcup_{R> 0} \mM^R_{(p,q,x);y}(f_1,f_2,g;H,J)
\end{align*}
which is a manifold of dimension ${\rm d}(p,q,x;y)= \mu(x)-\mu(y)+i_{f_1}(p)+i_{f_2}(q)-2n +1$ as the preimage of $\Delta\times\Delta := \{ (x,x;y,y)|\, (x,y)\in L^2 \}\subset (L\times L)^2$ by
\begin{align*}
{\rm Ev}\co \left\{
\hspace{-.2cm}\begin{array}{rl}
\mW^u_p(f_1,g)\times \mW^u_q(f_2,g)\times \mM_{x,y}(H,J) & \ra  (L\times L)^2\\
(\gamma_1,\gamma_2,u) &\mapsto  (\gamma_1(0),\gamma_2(0);\gamma_1(R),u(0,0))
\end{array}
\right.
\end{align*}
with the additional $R$ parameter. \\

Its $0$--dimensional component is compact, and we define
\begin{align*}
F(p,q,x) := \sum_{y|\,{\rm d}(p,q,x;y)=0} \#_2\mM_{(p,q,x);y}(f_1,f_2,g;H,J) \cdot y.
\end{align*}
Let $\mM^0_{(p,q,x);z}(f_1,f_2,g;H,J)$ be obtained from $\mM^R_{(p,q,x);z}(f_1,f_2,g;H,J)$ when $R$ converges to $0$. Notice that we get $\mM^+_{(p,q,x);z}(f_1,f_2,g;H,J)$, when $R$ converges to infinity. Hence, the boundary of the compactification of its $1$--dimensional component is the disjoint union of (see Figure \ref{modstr1})
\begin{align}
&\,\;\;\bigcup_{p'\in{\rm Crit}(f_1)} \;\;\mM_{p,p'}(f_1,g) \times \mM_{(p',q,x);y}(f_1,f_2,g;H,J)\label{modfloerbord1}\\
&\,\;\;\bigcup_{q'\in{\rm Crit}(f_2)} \;\;\mM_{q,q'}(f_2,g) \times \mM_{(p,q',x);y}(f_1,f_2,g;H,J)\label{modfloerbord2}\\
&\bigcup_{x'\in\mI(L,L';\eta,H)} \mM_{x,x'}(H,J) \times \mM_{(p,q,x');y}(f_1,f_2,g;H,J)\label{modfloerbord3}\\
&\bigcup_{y'\in\mI(L,L';\eta,H)} \mM_{(p,q,x);y'}(f_1,f_2,g;H,J) \times \mM_{y',y}(H,J)\label{modfloerbord4}\\
&\;\;\;\;\;\;\;\;\bigcup \; \mM^+_{(p,q,x);y}(f_1,f_2,g;H,J) \;\bigcup\; \mM^0_{(p,q,x);y}(f_1,f_2,g;H,J)\label{modfloerbord5}
\end{align}
Boundaries (\ref{modfloerbord1}), (\ref{modfloerbord2}), (\ref{modfloerbord3}) and (\ref{modfloerbord4}) (appearing when $R$ converges to a strictly positive real number), stand for $\del F+ F\del$. Hence $H_*(\varphi^+) = H_*(\varphi^0)$ and the first equality in (\ref{module}) holds.

\begin{figure}[h]
\centering
\psfrag{B1}[][][1]{\footnotesize(\ref{modfloerbord1})}
\psfrag{B2}[][][1]{\footnotesize(\ref{modfloerbord2})}
\psfrag{B3}[][][1]{\footnotesize(\ref{modfloerbord3})}
\psfrag{B4}[][][1]{\footnotesize(\ref{modfloerbord4})}
\psfrag{B5}[][][1]{\footnotesize(\ref{modfloerbord5})}
\psfrag{u}[][][1]{\footnotesize$u$}
\psfrag{p}[][][1]{\footnotesize$p$}
\psfrag{p'}[][][1]{\footnotesize$p'$}
\psfrag{q}[][][1]{\footnotesize$q$}
\psfrag{q'}[][][1]{\footnotesize$q'$}
\psfrag{y}[][][1]{\footnotesize$y$}
\psfrag{y'}[][][1]{\footnotesize$y'$}
\psfrag{x}[][][1]{\footnotesize$x$}
\psfrag{x'}[][][1]{\footnotesize$x'$}
\psfrag{r}[][][1]{\footnotesize$r$}
\psfrag{R}[][][1]{\footnotesize$R$}
\psfrag{R>infty}[][][1]{\footnotesize$R\ra +\infty$}
\psfrag{R>0}[][][1]{\footnotesize$R\ra 0$}
\psfrag{g1R}[][][1]{\footnotesize$\gamma^R_1$}
\psfrag{g2}[][][1]{\footnotesize$\gamma_2$}
\includegraphics{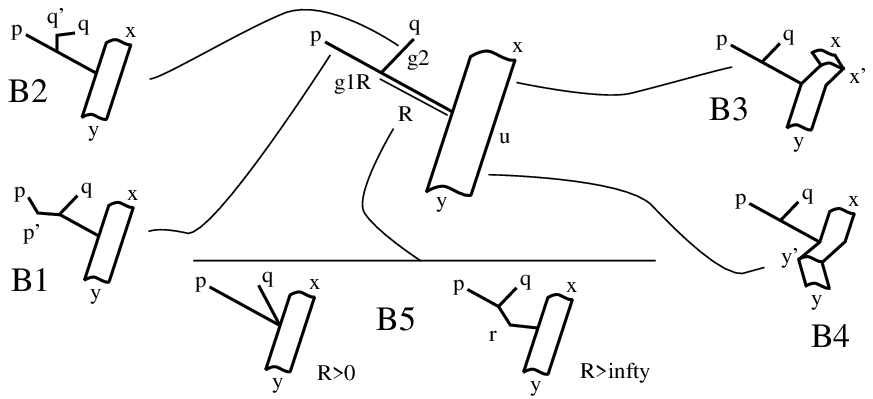}
\caption{Proof of the module structure (1)}
\label{modstr1}
\end{figure}

\hop {\it (step 2)} We construct other moduli spaces, similar to the previous ones:
\begin{align*}
\widehat{\mM}_{(p,q,x);y}(f_1,f_2,g;H,J):= \bigcup_{S< 0} \widehat{\mM}^S_{(p,q,x);y}(f_1,f_2,g;H,J)
\end{align*}
where $p\in {\rm Crit}(f_1)$, $q\in {\rm Crit}(f_2)$, $x$, $y\in\mI(L,L';\eta,H)$ and
\begin{align*}
&\widehat{\mM}^S_{(p,q,x);y}(f_1,f_2,g;H,J) := \\
& \left\{ (\gamma_1,\gamma^S_2,u)\in \mW^u_p(f_1,g)\times \mW^u_q(f_2,g)\times \mM_{x,y}(H,J) \left|
\begin{array}{l}
\gamma_1(0)=u(0,0) \\
\gamma^S_2(0)=u(S,0)
\end{array}
\right. \right\}.
\end{align*}
Their dimension is given by $\widehat{\rm d}(p,q,x;y)= i_{f_1}(p)+i_{f_2}(q) +\mu(x)-\mu(y)-2n+1$. By compactness, we define
\begin{align*}
G(p,q,x) := \sum_{y|\,\widehat{\rm d}(p,q,x;y)=0} \#_2 \widehat{\mM}_{(p,q,x);y}(f_1,f_2,g;H,J) \cdot y.
\end{align*}

The boundary of the (compactified) $1$--dimensional part, is the disjoint union
\begin{align}
&\,\;\;\bigcup_{p'\in{\rm Crit}(f_1)} \;\;\mM_{p,p'}(f_1,g) \times \widehat{\mM}_{(p',q,x);y}(f_1,f_2,g;H,J)\label{modfloerbord1'}\\
&\,\;\;\bigcup_{q'\in{\rm Crit}(f_2)} \;\;\mM_{q,q'}(f_2,g) \times \widehat{\mM}_{(p,q',x);y}(f_1,f_2,g;H,J)\label{modfloerbord2'}\\
&\bigcup_{x'\in\mI(L,L';\eta,H)} \mM_{x,x'}(H,J) \times \widehat{\mM}_{(p,q,x');y}(f_1,f_2,g;H,J)\label{modfloerbord3'}\\
&\bigcup_{y'\in\mI(L,L';\eta,H)} \widehat{\mM}_{(p,q,x);y'}(f_1,f_2,g;H,J) \times \mM_{y',y}(H,J)\label{modfloerbord4'}
\end{align}
when $S$ converges to a strictly negative real number, and
\begin{align}
& \bigcup \;\;\;\;\;\;\;\mM^-_{(p,q,x);y}(f_1,f_2,g;H,J)\label{modfloerbord5'}\\
& \bigcup \;\;\;\;\;\;\;\{ \Gamma\} \times \mM_{(p,x);y}(f_1,g;H,J)\label{modfloerbord6'}\\
& \bigcup \;\;\;\;\;\;\;\mM^0_{(p,q,x);y}(f_1,f_2,g;H,J)\label{modfloerbord7'}
\end{align}
when it converges respectively to $-\infty$ (\ref{modfloerbord5'}), (\ref{modfloerbord6'}) and $0$ (\ref{modfloerbord7'}) (see Figure \ref{modstr2}). Here, $\Gamma$ denotes the unique flow line of $f_2$ (if it exists), element of the unstable manifold of $q$ and reaching $x(0)$.

\begin{figure}[h]
\centering
\psfrag{B1}[][][1]{\footnotesize(\ref{modfloerbord1'})}
\psfrag{B2}[][][1]{\footnotesize(\ref{modfloerbord2'})}
\psfrag{B3}[][][1]{\footnotesize(\ref{modfloerbord3'})}
\psfrag{B4}[][][1]{\footnotesize(\ref{modfloerbord4'})}
\psfrag{B5}[][][1]{\footnotesize(\ref{modfloerbord5'})}
\psfrag{B6}[][][1]{\footnotesize(\ref{modfloerbord6'})}
\psfrag{B7}[][][1]{\footnotesize(\ref{modfloerbord7'})}
\psfrag{p}[][][1]{\footnotesize$p$}
\psfrag{p'}[][][1]{\footnotesize$p'$}
\psfrag{q}[][][1]{\footnotesize$q$}
\psfrag{q'}[][][1]{\footnotesize$q'$}
\psfrag{y}[][][1]{\footnotesize$y$}
\psfrag{y'}[][][1]{\footnotesize$y'$}
\psfrag{x}[][][1]{\footnotesize$x$}
\psfrag{x'}[][][1]{\footnotesize$x'$}
\psfrag{z}[][][1]{\footnotesize$z$}
\psfrag{S}[][][1]{\footnotesize$S$}
\psfrag{S>-infty}[][][1]{\footnotesize$S\ra -\infty$}
\psfrag{S>0}[][][1]{\footnotesize$S\ra 0$}
\psfrag{G}[][][1]{\footnotesize$\Gamma$}
\psfrag{g1}[][][1]{\footnotesize$\gamma_1$}
\psfrag{g2S}[][][1]{\footnotesize$\gamma_2^S$}
\psfrag{u}[][][1]{\footnotesize$u$}
\includegraphics{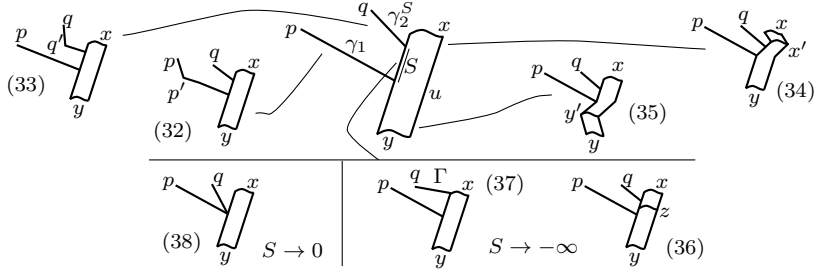}
\caption{Proof of the module structure (2)}
\label{modstr2}
\end{figure}

Boundaries (\ref{modfloerbord1'}), (\ref{modfloerbord2'}), (\ref{modfloerbord3'}) et (\ref{modfloerbord4'}) stand for $\del G+ G\del$. If $q=m$, the maximum of $f_2$, (\ref{rema:max}) allows us to conclude. Otherwise the unstable manifold of $q$ has codimension at least $1$. Hence, up to generic choices, the boundary (\ref{modfloerbord6'}) does not appear and $H_*(\varphi^0) = H_*(\varphi^-)$. The second equality of (\ref{module}) holds.
\end{proof}

\begin{rema}
{\it 1.} With standard methods, one can show that this construction does not depend on $f_1$, $f_2$, $g$, $H$ and $J$.\\
{\it 2.} It is well-known that semi-pairs of pants provide a product on Lagrangian Floer homology which gives an alternate definition of the same module structure:
\begin{align*}
\xymatrix{
HM_*(L)\otimes HF_*(L,L) \ar[r]^{\hspace{-.2cm}\phi_f^H\otimes \,\operatorname{id}} \ar[rd]_{\ast} & HF_*(L,L)\otimes HF_*(L,L) \ar[d]^{\rm Semi-PP} \\
& HF_*(L,L)
}
\end{align*}
However, semi-pairs of pants require distinct Hamiltonian perturbations at their ends. Thus the naturality morphism is not a priori compatible with the moduli spaces appearing in this description.
\end{rema}

\subsection{Naturality and Lagrangian PSS as morphisms of modules}\label{sec:independance}

In this section we prove that these two isomorphisms preserve the algebraic structures described above. More precisely, we prove the following lemma.
\begin{lemm}
For every $\alpha$ in $HM_*(L)$ and every $b$ in $HF_*(L,L')$, it holds
\begin{align}\label{mod0}
b_{H}(\alpha\ast b)= \alpha\ast b_{H}(b).
\end{align}
For every $\alpha$ and $\beta$ in $HM_*(L)$, it holds
\begin{align}\label{mod}
\phi_f^{H}(\alpha\cdot\beta)= \alpha\ast \phi_f^{H}(\beta).
\end{align}
\end{lemm}

\subsubsection{Naturality isomorphism}
Recall that, by naturality ($b_H^{-1}$)
\begin{itemize}
\item to $y$ and $z\in\mI(L,L';\eta,H)$, uniquely correspond $y'$ and $z'\in\mI(L,L'';\widetilde{\eta},0)$,
\item to $u\in \mM_{y,z}(L,L';H,J)$ uniquely corresponds $u'\in \mM_{y',z'}(L,L'';0,\widetilde{J})$.
\end{itemize}
Moreover $u'(0,0)=(\phi^0_H)^{-1}(u(0,0))=u(0,0)$. Hence we can deduce that $b_H^{-1}$ induces a one-to-one correspondence between moduli spaces (\ref{modpdt}) and transports the $HM_*(L)$--module structure of $HF_*(L,L';H,J)$, on the one of $HF_*(L,L'';0,\widetilde{J})$. Therefore (\ref{mod0}) is proved and $b_H$ is an isomorphism of $HM_*(L)$--modules.

\subsubsection{Lagrangian PSS isomorphism}

The proof of (\ref{mod}) is {\it formally} the same as the one of Lemma \ref{floermodule}. We choose two Morse functions $f_1$, $f_2$, a metric $g$ and a pair $(H,J)$ satisfying the same regularity conditions. First, note that the relation holds at the chain level if $\alpha$ is represented by the maximum $m$ of $f_1$.\\

\begin{figure}[h]
\centering
\psfrag{B1}[][][1]{\footnotesize(\ref{modpssbord1})}
\psfrag{B2}[][][1]{\footnotesize(\ref{modpssbord2})}
\psfrag{B3}[][][1]{\footnotesize(\ref{modpssbord3})}
\psfrag{B4}[][][1]{\footnotesize(\ref{modpssbord4})}
\psfrag{B5}[][][1]{\footnotesize(\ref{modpssbord5})}
\psfrag{p}[][][1]{\footnotesize$p$}
\psfrag{q}[][][1]{\footnotesize$q$}
\psfrag{y}[][][1]{\footnotesize$y$}
\psfrag{r}[][][1]{\footnotesize$r$}
\psfrag{R}[][][1]{\footnotesize$R$}
\psfrag{R>infty}[][][1]{\footnotesize$R\ra +\infty$}
\psfrag{R>0}[][][1]{\footnotesize$R\ra 0$}
\includegraphics{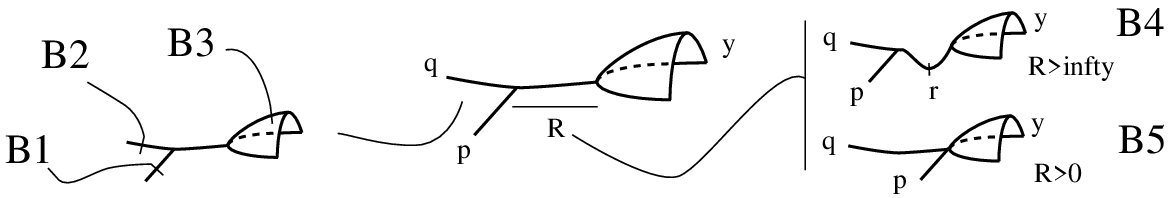}
\caption{Proof of structure preservation by the PSS morphism (1)}
\label{strucpreserv1}
\end{figure}

\hop {\it (Step $1$)} Now, $\mM_{(p,q);y}^R(f_1,f_2,g;H,J)$ is defined as the set
\begin{multline*}
\left\{ (\gamma_1,\gamma_2^R,u) \in \mW^u_p(f_1,g)\times \mW^u_q(f_2,g) \times \mM_y^H (J) \left|
\begin{array}{l}
\gamma_1(0)=\gamma_2^R(0) \\
\gamma_2^R(R)=u(-\infty,-)
\end{array} \right. \right\}.
\end{multline*}
We recall that $\mM_y^H (J)$ is a moduli space used in the definition of the PSS isomorphism (see \S \ref{subsec:pss} and Figure \ref{strucpreserv1}). The space $\mM^{0}_{(p,q);y}(f_1,f_2,g;H,J)$ is obtained by putting $R=0$ in the previous definition and
\begin{align*}
\mM^{+}_{(p,q);y}(f_1,f_2,g;H,J):=\bigcup_{r\in{\rm Crit}(f_1)} \;\mM_{(p,q);r}(f_1,f_2,f_2;g) \times \mM^{f_2,H}_{r,y}.
\end{align*}
Boundaries are now (see Figure \ref{strucpreserv1}):
\begin{align}
&\;\;\,\bigcup_{p'\in{\rm Crit}(f_1)} \;\;\mM_{p,p'}(f_1,g) \times \mM_{(p',q);y}(f_1,f_2,g;H,J)\label{modpssbord1}\\
&\;\;\,\bigcup_{q'\in{\rm Crit}(f_2)} \;\;\mM_{q,q'}(f_2,g) \times \mM_{(p,q');y}(f_1,f_2,g;H,J)\label{modpssbord2}\\
&\bigcup_{y'\in\mI(L,L';\eta,H)} \mM_{(p,q);y'}(f_1,f_2,g;H,J) \times \mM_{y',y}(H,J)\label{modpssbord3}\\
&\;\;\;\;\;\;\;\;\bigcup \;\;\;\;\;\;\;\mM^+_{(p,q);y}(f_1,f_2,g;H,J)\label{modpssbord4}\\
&\;\;\;\;\;\;\;\;\bigcup \;\;\;\;\;\;\;\mM^0_{(p,q);y}(f_1,f_2,g;H,J)\label{modpssbord5}
\end{align}
Boundaries (\ref{modpssbord1}), (\ref{modpssbord2}), and (\ref{modpssbord3}) appear when $R$ converges to a stricly positive real number and stand for $\del F+ F\del$, whereas (\ref{modpssbord4}) and respectively (\ref{modpssbord5}) appear when it converges to infinity and $0$. That shows $H_*(\varphi^+) = H_*(\varphi^0)$.\\

\hop {\it (Step $2$)} The set $\widehat{\mM}^S_{(p,q);y}(f_1,f_2,g;H,J)$ is defined as
\begin{multline*}
\left\{ (\gamma_1^S,\gamma_2,u)\in \mW^u_p(f_1,g)\times \mW^u_q(f_2,g)\times \mM_y^H (J) \left|
\begin{array}{l}
\gamma_1^S(0)=u(S,0) \\
\gamma_2(0)=u(-\infty,-)
\end{array}
\right. \right\}
\end{multline*}
and $\mM^{-}_{(p,q);y}(f_1,f_2,g;H,J):= \bigcup_{x\in\mI(L,L';\eta,H)} \,\mM_{q,x}^{f_2,H} \times \mM_{(p,x);y}(f_1,g;H,J)$. Notice that $\mM^{0}_{(p,q);y}(f_1,f_2,g;H,J)$ is now viewed as the limit of $\widehat{\mM}^S_{(p,q);y}(f_1,f_2,g;H,J)$ when $S$ converges to $-\infty$.

\begin{figure}[h]
\centering
\psfrag{B1}[][][1]{\footnotesize(\ref{modpssbord1'})}
\psfrag{B2}[][][1]{\footnotesize(\ref{modpssbord2'})}
\psfrag{B3}[][][1]{\footnotesize(\ref{modpssbord3'})}
\psfrag{B4}[][][1]{\footnotesize(\ref{modpssbord4'})}
\psfrag{B5}[][][1]{\footnotesize(\ref{modpssbord5'})}
\psfrag{B6}[][][1]{\footnotesize(\ref{modpssbord6'})}
\psfrag{p}[][][1]{\footnotesize$p$}
\psfrag{q}[][][1]{\footnotesize$q$}
\psfrag{y}[][][1]{\footnotesize$y$}
\psfrag{x}[][][1]{\footnotesize$x$}
\psfrag{S}[][][1]{\footnotesize$S$}
\psfrag{S>infty}[][][1]{\footnotesize$S\ra +\infty$}
\psfrag{S>0}[][][1]{\footnotesize$S\ra -\infty$}
\psfrag{G}[][][1]{\footnotesize$\Gamma$}
\includegraphics{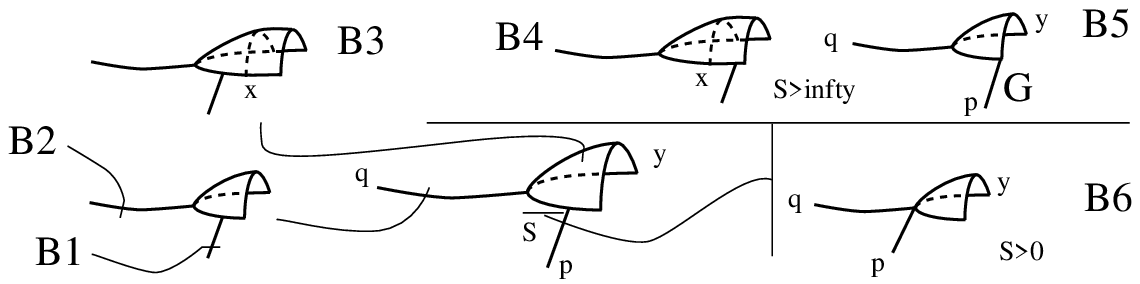}
\caption{Proof of structure preservation by the PSS morphism (2)}
\label{strucpreserv2}
\end{figure}

When $S$ converges to a real number, one gets (see Figure \ref{strucpreserv2})
\begin{align}
&\;\;\bigcup_{p'\in{\rm Crit}(f_1)} \;\;\mM_{p,p'}(f_1,g) \times \widehat{\mM}_{(p',q);y}(f_1,f_2,g;H,J)\label{modpssbord1'}\\
&\;\;\bigcup_{q'\in{\rm Crit}(f_2)} \;\;\mM_{q,q'}(f_2,g) \times \widehat{\mM}_{(p,q');y}(f_1,f_2,g;H,J)\label{modpssbord2'}\\
&\bigcup_{x\in\mI(L,L';\eta,H)} \,\widehat{\mM}_{(p,q);x}(f_1,f_2,g;H,J) \times \mM_{x,y}(H,J)\label{modpssbord3'}
\end{align}
which stand for $\del G+ G\del$. Otherwise, one gets
\begin{align}
&\bigcup \;\;\;\;\;\;\;\mM^{-}_{(p,q);y}(f_1,f_2,g;H,J)\label{modpssbord4'}\\
&\bigcup \;\;\;\;\;\;\;\{ \Gamma\} \times \mM_{q,y}^{f_2,H}\label{modpssbord5'}\\
&\bigcup \;\;\;\;\;\;\;\mM^0_{(p,q);y}(f_1,f_2,g;H,J)\label{modpssbord6'}
\end{align}
$S$ converging to $+\infty$ (\ref{modpssbord4'}), (\ref{modpssbord5'}) and to $-\infty$ (\ref{modpssbord6'}). Up to generic choices, $\Gamma$ (that is the unique flow line of $f_1$, element of the unstable manifold of $p$ reaching $y(0)$) does not exist. Hence (\ref{modpssbord5'}) does not appear and (\ref{mod}) is proved.

\subsubsection{Proof of Proposition \ref{prop:commut}}

Recall that we have to show that
\begin{align*}
\xymatrix{\relax
    HM_*(L;f,g) \ar[d]_{\phi_f^{H'}}\ar[r]^{\hspace{-.3cm}\phi_f^{H}} & HF_*(L,L;H,J) \ar[d]^{b_H^{-1}}\\
    HF_*(L,L;H',J') \ar[r]_{b_{H'}^{-1}} & HF_*(L,L_0;0,\widetilde{J})
}
\end{align*}
is commutative, as soon as $\phi_H^*J=\phi_{H'}^*J'=:\widetilde{J}$ and $(\phi_H^1)^{-1}(L)=(\phi^1_{H'})^{-1}(L)=:L_0$.

\begin{proof}[Proof of Proposition \ref{prop:commut}]
For such two pairs $(H,J)$ and $(H',J')$, we define
\begin{align*}
\Phi := (\phi^{H'}_f)^{-1} \circ b_{H'} \circ b_H^{-1} \circ \phi_f^H \co HM_*(L)\longrightarrow HM_*(L)
\end{align*}
which is an isomorphism. Since $[L]$ generates $HM_n(L)$, $\Phi([L])= [L]$ (due to $\bZ_2$ coefficients, signs are arbitrary). Moreover, $[L]$ is the unit of the {\it ring} $(HM_*(L),\cdot)$ and $\Phi$ preserves the {\it module} structure of $HM_*(L)$. Hence we have
\begin{align*}
\Phi(a) &=\Phi(a\cdot [L]) = a\cdot \Phi([L]) =a\cdot [L] =a.
\end{align*}
Therefore, $\Phi$ is the identity and the proposition is proved.
\end{proof}

\begin{rema}\label{remarque:reference action et maslov diagramme commutativif}
{\it 1.} The PSS morphism does not require any particular choice of the action reference $\eta$. On the other hand, our choices have to be consistent with the naturality morphism. Hence we choose an element $\eta_0\in L\cap L_0$ and we put $\eta := b_H(\eta_0)$ and $\eta' := b_{H'}(\eta_0)$ as action references for $HF(L,L;H,J)$ and $HF(L,L;H',J')$.

{\it 2.} The PSS morphism requires a particular choice of references for the Maslov index, $z_0\in \mI(L,L;\eta,H)$ and $z'_0\in \mI(L,L;\eta',H')$. If we denote by $\widetilde{z}_0$ and $\widetilde{z}'_0$ the elements of $\mI(L,L_0;\eta_0,0)$ such that $b_{H}(\widetilde{z}_0)=z_0$ and $b_{H'}(\widetilde{z}'_0)=z'_0$, $\Phi$ is more precisely a map from $HM_*(L)$ to $HM_{*+\delta}(L)$ where $\delta:=\mu(\widetilde{z}_0,\widetilde{z}'_0)$. But $\Phi$ is injective and thus $\Phi([L])$ and $\Phi(1)$ have to be non zero. Hence $\delta=0$ and $\Phi$ preserves the grading.
\end{rema}

\subsection{End of the proof of Theorem \ref{theo:independ}}\label{sec:continuity}

We first give a standard result (Lemma \ref{lemm:estime} below) used in the Hamiltonian case by Schwarz \cite{Schwarz00}. We recall that a time-dependent Hamiltonian function $H$ is said to be normalized if it satisfies
\begin{align*}
\forall t\in I, \;\; \int_M H(t,-)\,\omega^n=0
\end{align*}
if $M$ is compact, or if it has compact support otherwise. The quantities
\begin{align*}
E_-(H):=\int_I \min_M H \,dt \;\;{\rm and}\;\; E_+(H):=\int_I \max_M H \,dt
\end{align*}
give rise to Hofer's norm on the set of normalized Hamiltonian functions and Hofer's distance on the set of Lagrangian submanifolds Hamiltonian isotopic to a fixed one:
\begin{align*}
\| H \|:= E_+(H) - E_-(H)  \;\;{\rm and}\;\; \nabla(L,L'):= \inf \{ \| H \| |\, H\;{\rm satisfies}\; \phi^1_H(L) =L' \}.
\end{align*}
We will need the following lemma, in order to complete the proof of Theorem \ref{theo:independ}.
\begin{lemm}\label{lemm:estime}
For regular pairs $(H,J)$, $(H',J')$ and $0\neq\alpha\in H_*(L)$, we have
\begin{align}\label{estime1}
E_-(H'-H) + \alpha_{\eta,\eta'} \leq \sigma_L(\alpha;H',J',\eta') - \sigma_L(\alpha;H,J,\eta) \leq E_+(H'-H) + \alpha_{\eta,\eta'},
\end{align}
where $\alpha_{\eta,\eta'}=\int u^* \omega$ for any map $u\co I\times I\ra M$ such that $u(I,\{ 0,1 \})\subset L$, $u(0,-)=\eta$ and $u(1,-)=\eta'$. Moreover, the following inequality holds
\begin{align}\label{estime}
|c_L(\alpha;H',J') - c_L(\alpha;H,J)| \leq \|H'-H \|.
\end{align}
\end{lemm}
For clarity, its proof will be given at the end of this section (\S \ref{subsec:lemm}). In the next two paragraphs, we combine Proposition \ref{prop:commut} and Lemma \ref{lemm:estime} to prove the invariance part (\S \ref{subsec:inv}) and the continuity part (\S \ref{subsec:cont}) of Theorem \ref{theo:independ}.

\subsubsection{Proof of the invariance part of Theorem \ref{theo:independ}}\label{subsec:inv}
Since the naturality morphism preserves the action for our choices of references (see Remark \ref{remarque:reference action et maslov diagramme commutativif}), Proposition \ref{prop:commut} implies
\begin{align*}
\sigma_L(-;H,J,\eta) &= \widetilde{\sigma}_{L,L_0}(b_H^{-1} \circ \phi^H_f(-);0,\widetilde{J},\eta_0)\\
 &= \widetilde{\sigma}_{L,L_0}(b_{H'}^{-1} \circ \phi^{H'}_f(-);0,\widetilde{J},\eta_0) = \sigma_L(-;H',J',\eta')
\end{align*}
whenever $\phi_H^*J=\phi_{H'}^*J'=:\widetilde{J}$ and $(\phi_H^1)^{-1}(L)=(\phi^1_{H'})^{-1}(L)=:L_0$. Moreover, as an obvious corollary of Lemma \ref{lemm:estime}, putting $H=H'$ and $\eta=\eta'$ in (\ref{estime1}), we get that the homological Lagrangian spectral numbers (relative and absolute versions) do not depend on the almost complex structure. Hence their absolute version only depend on $L$ and $L_0$. We put
\begin{align*}
c(-;L,L_0):= c_L(-;H,J)= \sigma_L(-;H,J,\eta) - \sigma_L(1;H,J,\eta)
\end{align*}
for any Hamiltonian function $H$ such that $\phi_H^1(L_0)=L$.

\begin{rema}
Notice that, as mentioned by McDuff and Salamon \cite{McDuffSalamon04}, this method can be used in the Hamiltonian case. Indeed, in that case, the module structure is known (introduced by Floer \cite{Floer89} and studied by Le and Ono \cite{LeOno96}), and the naturality morphism is a particular case of isomorphisms described by Seidel \cite{Seidel97}.
\end{rema}

\subsubsection{Proof of the continuity part of Theorem \ref{theo:independ}}\label{subsec:cont}
This property is a corollary of Lemma \ref{lemm:estime}.
\begin{coro}
Let $L$, $L_0$ and $L_1$ be Hamiltonian isotopic Lagrangian submanifolds of $M$. We have, for any homology class of $L$, $\alpha\neq 0$:
\begin{align}\label{cont}
| c(\alpha;L,L_0) - c(\alpha;L,L_1) | & \leq\nabla (L_0,L_1)
\end{align}
which in particular implies that $c(\alpha;L,-)$ is continuous on the set of the Lagrangian submanifolds Hamiltonian isotopic to $L$, with respect to Hofer's distance.
\end{coro}
Indeed, let $G$ and $H$ be normalized Hamiltonian functions such that $\phi^1_{G}(L_1)=L_0$ and $\phi^1_H(L_0)=L$. Then the Hamiltonian diffeomorphism $\psi:= \phi^1_H \circ \phi^1_G$ satisfies $\psi^{-1}(L)=L_1$. Moreover it is the end of the Hamiltonian isotopy induced by:
\begin{align*}
H'_t := H_t + G_t \circ (\phi^t_H)^{-1}.
\end{align*}
Since $\phi_H^t$ is a symplectomorphism, Lemma \ref{lemm:estime} implies that:
\begin{align*}
| c(\alpha;L,L_0) - c(\alpha;L,L_1) | & \leq \| H- H'\|
=  \| G \circ (\phi_H)^{-1} \|
=  \| G \|.
\end{align*}
As the left side does not depend on $G$ (as long as $\phi^1_{G}(L_1)=L_0$), we get (\ref{cont}).

\subsubsection{Proof of Lemma \ref{lemm:estime}}\label{subsec:lemm}
It remains to prove this lemma in order to end the proof of Theorem \ref{theo:independ}
\begin{proof}
Let $(H,J)$ and $(H',J')$ be two regular pairs. We consider the classical comparison morphism induced by the homotopy $\widetilde{G}_s:=H + \beta(s)(H'-H)$, where $\beta$ is a smooth, increasing function whose value is $0$ for $s<-1$ and $1$ for $s>1$. If $\psi^{H,H'}(x)=\sum_{y} a_{y}\cdot y$ with $a_y\neq 0$ for an element $x\in \mI(L,L;\eta,H)$ then there exists a Floer trajectory $u$ connecting $x$ and $y$. Its energy is easily seen to satisfy
\begin{align*}
\int_{{\Bbb R}\times I} \| \del_s u\|^2 = \mA_H(x) - \mA_{H'}(y) + \int_{{\Bbb R}} \beta'(s) \left( \int_I (H_t'-H_t)\circ u(s,t)\, dt\right) ds + \alpha_{\eta,\eta'}
\end{align*}
and hence
\begin{align}\label{estcomp}
\mA_{H'}(y)\leq \mA_H(x) + E_+(H'-H) + \alpha_{\eta,\eta'}.
\end{align}

Therefore for any non zero element $x$ in $HF_*(L,L;H,J)$ we get
\begin{align*}
\widetilde{\sigma}_{L,L}(\psi^{H,H'}(x);H',J',\eta') \leq \widetilde{\sigma}_{L,L}(x;H,J,\eta) + E_+(H'-H) + \alpha_{\eta,\eta'}.
\end{align*}
Lemma \ref{lemm:PSScomparaison} implies that $\phi^{H'}_f = \psi^{H,H'}\circ \phi^H_f$. Thus for $\alpha\neq 0$ in $HM_*(L)$, we obtain
\begin{align*}
\widetilde{\sigma}_{L,L}(\phi_f^{H'}(\alpha);H',J',\eta') &= \widetilde{\sigma}_{L,L}(\psi^{H,H'}\circ \phi^H_f(\alpha);H',J',\eta') \\
& \leq \widetilde{\sigma}_{L,L}(\phi^H_f(\alpha);H,J,\eta) + E_+(H'-H) + \alpha_{\eta,\eta'}.
\end{align*}
Hence we have $\sigma_L(\alpha;H',J',\eta') \leq \sigma_L(\alpha;H,J,\eta) +E_+(H'-H) +\alpha_{\eta,\eta'}$. Inverting $H$ and $H'$, $\eta$ and $\eta'$ and using the equalities $E_+(-K)=-E_-(K)$ and $\alpha_{\eta,\eta'}=-\alpha_{\eta',\eta}$ give (\ref{estime1}). Writing (\ref{estime1}) for $\alpha=1$ and subtracting it to (\ref{estime1}) give (\ref{estime}).
\end{proof}

\subsection{Recovering Hamiltonian spectral invariants}\label{sec:classical}
We show here that the Hamiltonian spectral invariants can be viewed as particular Lagrangian ones. In order to do so, we prove equality (\ref{prop:generalisation}) which states that they coincide via the Biran--Polterovich--Salamon isomorphism \cite{BiranPolterovichSalamon03}. We first recall the construction of this morphism which compares the Hamiltonian Floer homology of a symplectic manifold $(M,\omega)$ and the Lagrangian Floer homology of $(\underline{M}, \underline{\omega}):=(M\times M,\omega\oplus (-\omega))$ with respect to its Lagrangian submanifold $\Delta:= \{ (x,x) \in M\times M \,|\; x\in M\}$.\\

Let $(H,J)$ be a regular pair on $(M,\omega)$. We define $\underline{H}$ and $\underline{J}$ on $\underline{M}$ by
\begin{align*}
&\underline{H}_t(x_1;x_2) := H_t(x_1) + H_{1-t}(x_2)\;\;{\rm and}\\
&\underline{J}_t := (\underline{\phi^t})^*(J_t\times -J_{1-t})= ((\phi_H^t)^* J_t \times -(\phi^{1-t}_H\circ (\phi^1_H)^{-1})^* J_{1-t})
\end{align*}
with $t\in [0, 1/2]$. In order to shorten notation, we denote Hamiltonian vector fields and flows by $X_t$ and $\phi^t$. Those induced by $\underline{H}$ are given by
\begin{align*}
\underline{X}_t(x_1,x_2) &= (X_t(x_1),-X_{1-t}(x_2)),\\
\underline{\phi}^t(x_1,x_2) &= (\phi^t(x_1),\phi^{1-t}\circ (\phi^1)^{-1}(x_2)).
\end{align*}

Since $\phi^1$ is a symplectomorphism, its graph, denoted by $\Gamma_{\phi^1}$, is a Lagrangian submanifold of $M\times M$. The generators of the Floer complex of $(M,\omega)$, $\mP^0(H)$, and of the Lagrangian Floer complex of $(\underline{M},\underline{\omega})$, with respect to $\Delta$ and $\Gamma_{\phi^1}$, are identified by the bijection
\begin{align*}
\mP^0(H) \ni [ x(t)\co t\mapsto \phi^t(x(0)) ] \leftrightarrow \underline{x}:=(x(0),x(0)) \in \Delta\cap \Gamma_{\phi^1}.
\end{align*}
Moreover, the transversality conditions required in the Lagrangian case ($\Delta\cap\Gamma_{\phi^1}$ is transverse) and in the Hamiltonian case ($\det (d\phi^1 - {\rm id})\neq 0$) are equivalent.

To $u\co\bR\times S^1\ra M$, let us associate $\underline{u}\co\bR\times [0,1/2]\ra \underline{M}$, defined by
\begin{align*}
\underline{u}(s,t):= (\underline{\phi^t})^{-1}(u(s,t),u(s,1-t))=((\phi^t)^{-1}(u(s,t)), \phi^1\circ(\phi^{1-t})^{-1}(u(s,1-t))).
\end{align*}
Notice that we have $\underline{u}(s,0)\in \Delta$ and $\underline{u}(s,1/2) \in \Gamma_{\phi^1}$. Projecting the relation $\del_s \underline{u} +\underline{J}\del_t \underline{u} =0$ gives the following equivalences:
\begin{align*}
\del_s u +J\del_t u +\nabla H_t(u)=0 &\;\;\;\mbox{iff}\;\;\; \del_s \underline{u} +\underline{J}\del_t \underline{u} =0, \\
u(\pm \infty,t)=x^\pm (t),\; E(u)<\infty &\;\;\;\mbox{iff}\;\;\; \underline{u}(\pm \infty,t)=\underline{x}^\pm,\; E(\underline{u})<\infty.
\end{align*}
Since $\del_s \underline{u} (s,t) = d\underline{\phi}^t (\del_s u(s,t), \del_s u(s,1-t))$, the energy is preserved, in the sense that $E(\underline{u})= E(u)$. Via (\ref{actionenergie}), the action is preserved up to an additive constant. Hence we get an identification between the two chain complexes:
\begin{align*}
(CF_*(M,\omega;H,J),\del) \simeq (CF_*(\underline{M},\underline{\omega};\Delta,\Gamma_{\phi^1};0,\underline{J}),\underline{\del})
\end{align*}
which induces the desired isomorphism in homology.\\

Let $(f,g)$ be a Morse--Smale pair on $M$. The pair $(\underline{f},\underline{g})$, defined on $\Delta\subset M\times M$ by $\underline{f}(x,x):=f(x)$ and $\underline{g}_{(x,x)}((v,v),(w,w)):=g_{x}(v,w)$, is Morse--Smale. Since the critical points of $f$ and its flow lines with respect to $g$ are identified with the critical points of $\underline{f}$ and its flow lines with respect to $\underline{g}$, the above process leads to an identification between the moduli spaces defining the PSS morphisms in both cases. Therefore the following diagram is commutative:
\begin{align*}
\xymatrix{\relax
    HM_*(M;f,g) \ar@{=}[r] \ar[d]_{\rm PSS} & HM_*(\underline{M};\underline{f},\underline{g}) \ar[d]^{{\rm naturality}\;\circ\; {\rm PSS}}\\
    HF_*(M,\omega;H,J) \ar@{=}[r]^{\hspace{-.5cm}{\rm BPS}} & HF_*(\underline{M},\underline{\omega};\Delta,\Gamma_{\phi^1};0,\underline{J})
}
\end{align*}
Together with the fact that the action is preserved up to a constant, this gives the correspondence between spectral invariants in both cases. The equality (\ref{prop:generalisation}) $\rho(\alpha;\phi) - \rho(1;\phi) = c(\underline{\alpha};\Delta,\Gamma_{\phi})$ holds as claimed in the introduction.


\section{Lagrangian spectral invariants of higher order}\label{sec:3}

In this section, we use the Barraud--Cornea spectral sequence machinery to introduce Lagrangian spectral invariants of higher order, which generalize the invariants defined before. In order to do so, we first recall this machinery. Then we define the generalized invariants and prove their main property. Details and proofs about this spectral sequence machinery can be found in \cite{BarraudCornea05, BarraudCornea07}.


\subsection{Recollection of the Barraud--Cornea spectral sequences}

At the center of this construction, there is a chain complex whose differential takes into account higher dimensional moduli spaces of Morse (or Floer) trajectories (not only their $0$--dimensional component). The relevant spectral sequence is induced by a natural filtration of this extended chain complex.

\subsubsection{Morse case}

Let $L^n$ be a smooth manifold and $f$ a Morse function. Let $w$ be a path embedded in $L$, which passes through all the critical points of $f$. We denote by $\rho$ the projection of $L$ onto $\widetilde{L}:=L/ w$ (the critical points of $f$ are mapped to a distinguished point, $*$, of $\widetilde{L}$). Let $X$ be simply-connected and $l\co\widetilde{L}\ra X$ be a continuous map. Notice that $L$ and the quotient space $\widetilde{L}$ are clearly homotopy equivalent via $\rho$. The homotopy class of the map $l$ (rather than the map itself) being involved in the following construction, one can actually view $l$ as a continuous map from $L$ to $X$.\\

Let us denote by $C_{p,q}L$ the subset of $C^0([0,f(p)-f(q)],L)$ consisting of paths from $p$ to $q$ and by $\gamma_{p,q}$ the map from $\mM_{p,q}(f,g)$ to $C_{p,q}L$ which associates to $x$ in the connecting manifold of $p$ and $q$, the flow line representing it, parameterized by $[0,f(p)-f(q)]$. As this map is compatible with the compactification of moduli spaces, it induces $\overline{\gamma}_{p,q}\co \overline{\mM}_{p,q}(f,g)\ra C_{p,q}L$.

Let $\Omega'Y$ be the space of Moore loops of $Y$, that is, the set of loops in $Y$ parameterized by any compact interval $[0,a]$, with $a>0$. We denote by $\mS_*(\Omega'Y)$ the cubical chain complex of $\Omega'Y$. Concatenation in $\Omega'Y$ defines a product on $\mS_*(\Omega'Y)$. We define $\mQ_{p,q}\co C_{p,q}L\ra \Omega' \widetilde{L}$ as the projection by $\rho$ of the paths in $L$. We denote by $\Phi_{p,q}$ the composition $\mQ_{p,q}\circ\overline{\gamma}_{p,q}$. The map $l$ induces maps
\begin{align*}
\Omega'l\co \Omega'\widetilde{L}\ra \Omega'X\;\;\mbox{ and }\;\;(\Omega'l)_*\co \mS_*(\Omega' \widetilde{L}) \longrightarrow \mS_*(\Omega'X)
\end{align*}
preserving the multiplication. Notice that $\mR_*^X:= \mS_*(\Omega'X)$ is a differential ring.

We denote by $j_*$ the map induced at the (cubical) chain level by the inclusion of $(\overline{\mM}_{p,q}(f,g),\emptyset)$ in $(\overline{\mM}_{p,q}(f,g),\del \overline{\mM}_{p,q}(f,g))$. The last ingredient is a representing chain system of $\overline{\mM}(f,g)$, that is, a set
\begin{align*}
\zeta=\{ s_{pq}\in \mS_{i_f(p)-i_f(q)-1}(\overline{\mM}_{p,q}(f,g)) \, | \; p,\, q\in {\rm Crit}(f)\}
\end{align*}
such that
\begin{itemize}
\item[i.] the chain $\sigma_{pq}:=j_*(s_{pq}) \in \mS_{i_f(p)-i_f(q)-1}(\overline{\mM}_{p,q}(f,g),\del \overline{\mM}_{p,q}(f,g))$ represents the fundamental class of $\overline{\mM}_{p,q}(f,g)$ relative to the boundary,
\item[ii.] $\del s_{pq}= \sum_r s_{pr}\times s_{rq}$.
\end{itemize}

\hop Such a system exists for $\overline{\mM}(f,g)$. It can be built by induction, using ii.\\

We denote by $a_{pq}:=(\Omega'l)_*\circ(\Phi_{p,q})_*(s_{pq})\in \mR_{*}^X$. We form the left $\mR_*^X$--module $\mC_*^X=\mR_*^X\otimes CM_*(f,g)$ and we define a differential by the formula
\begin{align*}
\del(p)=\sum_{q\in {\rm Crit}f} a_{pq}\otimes q.
\end{align*}
Let $d'$ be the unique differential on $\mC_*^X$, satisfying $d'(a\otimes p)=(da)\otimes p+ a\cdot(\del p)$. Barraud and Cornea showed that $d'^2=0$. Thus $(\mC_*^X,d')$ is a graded differential module. It is called the extended Morse complex and admits the following filtration:
\begin{align*}
F^k\mC^X= \mR_*^X\otimes \langle {\rm Crit}_j(f) \, | \; j\leq k\rangle_{\bZ_2} = \bigoplus_{j\leq k} \mR_*^X\otimes CM_j(f,g).
\end{align*}
The spectral sequence which we are interested in is the one based on this filtration:
\begin{align*}
EM(L;f,g;X) = (EM^r_{p,q}(L;f,g;X),d^r).
\end{align*}

Recall from the introduction that $E_X(L)$ denotes the Serre spectral sequence associated to the fibration $\Omega X\hookrightarrow E\ra L$ (pull-back of the path-loop fibration of $X$ over $l$). It is shown in \cite{BarraudCornea07} that there is an isomorphism of spectral sequences between $EM(L;f,g;X)$ and $E_X(L)$ at page $2$.

\subsubsection{Lagrangian Floer case}
This construction is analogous to the previous case. $w$ is an embedded path in $L$ which passes through $x(0)$ for every $x\in\mI(L,L';\eta,H)$. We denote by $\widetilde{M}:=M/{\rm Im}(w)$, by $\widetilde{L}:=L/{\rm Im}(w)$ and by $\rho$ the projections:
\begin{align*}
M\ra\widetilde{M},\;\; L\ra\widetilde{L},\;\; \mbox{ and } \;\; \mP_\eta (L,L') \ra \mP_\eta (\widetilde{L},L').
\end{align*}
Now $C_{x,y}\mP$ stands for the subset of $C^0([0,\mA_H(x)-\mA_H(y)],\mP_\eta(L,L'))$ consisting of paths with respective endpoints $x$ and $y$ and a map
\begin{align*}
\bar{\gamma}_{x,y}\co \overline{\mM}_{x,y}(L,L';H,J)\ra C_{x,y}\mP
\end{align*}
is defined. The evaluation map at $0$ is used to define $\mQ_{x,y}={\rm ev}_0 \circ \rho \co C_{x,y}\mP \ra\Omega'\widetilde{L}$ and $\Phi_{x,y}\co \overline{\mM}(x,y) \ra \Omega'\widetilde{L}$ denotes $\mQ_{x,y}\circ\bar{\gamma}_{x,y}$. Critical points being replaced by orbits and Morse index by Maslov index, the end of the procedure is formally identical and leads to the extended Floer complex $(\mC_*^X, d')$ which admits a filtration
\begin{align*}
F^k\mC^X= \mR_*^X\otimes \langle x\in \mI(L,L';\eta,H) \, | \;\mu(x)\leq k\rangle_{\bZ_2} = \bigoplus_{j\leq k} \mR_*^X\otimes CF_j (L,L';H,J).
\end{align*}
We deduce the desired spectral sequence:
\begin{align*}
EF(L,L';\eta;H,J;X) = (EF^r_{p,q}(L,L';\eta;H,J;X),d^r).
\end{align*}

Two of its properties (see \cite{BarraudCornea07}) will be useful in the next section:
\begin{enumerate}
\item $EF^2_{p,q}(L,L';\eta;X)\simeq H_q(\Omega X)\otimes HF_p(L,L';\eta)$ for any integers $p$ and $q$,
\item if $d^r \neq 0$, there exist $x$ and $y\in \mI(L,L';\eta,H)$ with $\mu(x,y)\leq r$ such that $\mM_{x,y}(L,L';H,J)$ is not empty.
\end{enumerate}

The comparison and naturality morphisms between (suitable) Floer complexes naturally extend to the extended complexes and induce isomorphisms, $\Psi^{H,H'}$ and ${\rm B}_H$, between the respective spectral sequences. Moreover, when $L'=L$, there also exists an extension of the Lagrangian PSS morphism to the Morse and Floer extended complexes, inducing an isomorphism $\Phi_f^H$ between the respective spectral sequences at page $2$. It restricts itself to the PSS morphism on $EM^2_{*,0}(L;f,g;X)\simeq HM_*(L;f,g)$. For these three constructions (see~\cite{BarraudCornea05}), one has to include the components of any dimension of the moduli spaces (not only the $0$--dimensional part). Therefore, in particular, Morse and Floer spectral sequences are identified at page $r$ ($r\geq 2$) with $E_X(L)$.

\subsection{Definition of Lagrangian spectral invariants of higher order}

For $\nu\in\bR$, we define a truncated, extended Floer complex by
\begin{align*}
\mC^{\nu} = \mC^{\nu}(L,L;H,J) := \mR^X_*\otimes CF_*^\nu(L,L;H,J).
\end{align*}
We denote by $\mC^{>\nu}$ the quotient $\mC / \mC^\nu$. Since the filtration is preserved, they give rise to spectral sequences denoted by $EF_{\nu}(L,L;H,J;X)$ and $EF^{>\nu}(L,L;H,J;X)$. Moreover, the inclusion and quotient which form the following short exact sequence
\begin{align*}
\xymatrix{
0 \ar[r] & \mC^\nu \ar@{^{(}->}[r]^i & \mC \ar@{>>}[r]^{\hspace{-.3cm}q} & \mC / \mC^\nu \ar[r] & 0
}
\end{align*}
are morphisms of filtered graded differential modules. Hence they induce morphisms between the respective spectral sequences
\begin{align}\label{eq:suite}
\xymatrix{
EF_{\nu}(L,L;H,J;X) \ar[r]^{i_\nu} & EF(L,L;H,J;X) \ar[r]^{\hspace{-.2cm}q_\nu} & EF^{>\nu}(L,L;H,J;X)
}.
\end{align}
\begin{defi}\label{defi:si2}
Let $\alpha\neq 0$ be an element of $EM^r_{*,*}(L;f,g;X)$. Its associated relative higher order spectral numbers are
\begin{align*}
\sic{r}{\alpha} &:= \inf \{ \nu \in\bR |\, \Phi_f^H(\alpha) \in {\rm Im}(i_\nu) \},\\
\skc{r}{\alpha} &:= \inf \{ \nu \in\bR |\, q_\nu(\Phi_f^H(\alpha))=0 \}.
\end{align*}
\end{defi}
It is useful to notice that (\ref{eq:suite}) is not a priori exact except for the second page (where the vector spaces appearing in the spectral sequences are essentially homological). However, we obviously still have $q_\nu\circ i_\nu=0$. Hence, in general
\begin{align*}
\skc{r}{-} \;\leq \;\sic{r}{-}
\end{align*}
with equality (at least) for $r=2$. In that case, let us choose a basis, $\{x_i\}$, of $H_q(\Omega X)$ and let $\alpha_j$ be non zero elements in $H_p(L)$. We put $\alpha:=\sum_j x_j\otimes \alpha_j$ in $EM^2_{p,q}(L;f,g;X)$. We have
\begin{align}\label{eq:egalite teutone}
\sic{2}{\alpha}=\skc{2}{\alpha}= \max_j \{ \sigma_L(\alpha_j;H,J,\eta) \}.
\end{align}
In particular, if $1$ denotes the generator of $H_0(\Omega X)$, for $\beta\neq 0$ in $H_*(L)$, we have $\sic{2}{1\otimes\beta}=\skc{2}{1\otimes\beta}= \sigma_L(\beta;H,J,\eta)$.\\

Like in the homological case, we now define an absolute version which does not depend on the choice of the reference path of the action. Since $EM(L;f,g;X)$ is a first quadrant spectral sequence, for all integer $r$, we have
\begin{align*}
EM^r_{0,0}(L;f,g;X) \simeq EM^2_{0,0}(L;f,g;X) \simeq H_0(\Omega X)\otimes H_0(L).
\end{align*}
Hence, it contains an element, denoted $1_r$, which corresponds to $1\otimes 1$ by the previous isomorphism. We use this element in order to normalize $\sic{r}{-}$ and $\skc{r}{-}$, using (\ref{eq:egalite teutone}) with $\alpha=1\otimes 1$.
\begin{defi}\label{defi:si2bis}
Let $\alpha\neq 0$ be an element of $EM^r_{*,*}(L;f,g;X)$. Its associated (absolute) higher order spectral invariants are
\begin{align*}
\ic{r}{\alpha} &:= \sic{r}{\alpha} - \sigma_L(1;H,J,\eta),\\
\kc{r}{\alpha} &:= \skc{r}{\alpha} - \sigma_L(1;H,J,\eta)
\end{align*}
where $1$ is the generator of $H_0(L)$.
\end{defi}

Notice that this absolute version still satisfies equalities (\ref{eq:egalite teutone}).

\begin{rema}
As suggested by our notation, these numbers satisfy the same invariance property as their homological counterpart. Indeed, the commutative diagram of Proposition \ref{prop:commut} induces a commutative diagram for the respective spectral sequences and the estimate (\ref{estcomp}) (and thus Lemma \ref{lemm:estime}) holds for them.
\end{rema}

\subsection{Spectral invariants of higher order of $(S^2\times S^4)\# (S^2\times S^4)$}

We compute explicitly the higher order (Morse) spectral invariants of $(S^2\times S^4)\# (S^2\times S^4)$. First we consider $S^2\times S^4$ with the Morse function $f$ defined as the sum of the height function on each sphere:
\begin{align*}
f_n \co S^n\lra \bR &\;\mbox{ such that }\; f(x_1,\ldots x_n)=x_n,\\
f \co S^2\times S^4\lra \bR &\;\mbox{ such that }\; f(x,y)= f_2(x) + f_4(y).
\end{align*}
The function $f$ has $4$ critical points, $p_6$, $p_4$, $p_2$, $p_0$ (namely $(\max(f_2),\max(f_4))$, $(\min(f_2),\max(f_4))$, $(\max(f_2),\min(f_4))$ and $(\min(f_2),\min(f_4))$) where $i_f (p_i)=i$.

The non empty moduli spaces of flow lines are the following
\begin{align*}
\mM_{p_6,p_4}, \; \mM_{p_6,p_2}, \; \mM_{p_6,p_0} \;\; {\rm and} \;\; \mM_{p_4,p_0}, \; \mM_{p_2,p_0}.
\end{align*}

We denote by $\alpha$ and $\beta$ the generators of $\mS_1(\Omega S^2)$ and $\mS_3(\Omega S^4)$ (seen as elements of $\mS_*(\Omega'(S^2\times S^4))$) and by $\gamma$ their product. They induce a representing chain system for the moduli spaces of Morse flow lines such that the (complete) differential of the extended Morse complex gives
\begin{align*}
\del p_6 = \alpha\otimes p_4 + \beta \otimes p_2 +\gamma\otimes p_0 \;\; {\rm and} \;\; \del p_4 = \beta \otimes p_0, \; \del p_2 = \alpha \otimes p_0, \; \del p_0 = 0.
\end{align*}
The element $\alpha\otimes p_4$ is an element of the second page of the Morse spectral sequence and by definition of the higher order (Morse) spectral invariants, we have
\begin{align}\label{c2}
\overline{c}^2(\alpha\otimes p_4)= c(p_4) = f(p_4)
\end{align}
(notice that $p_4$ denotes here the critical point and the induced homology class).\\

Let us now consider $(S^2\times S^4)_{[1]}\# (S^2\times S^4)_{[2]}$ with the Morse function respecting the two components \textquotedblleft$f= f_{[1]} \# f_{[2]}$" and extended on the connected sum (the notation $[i]$ is used in order to indicate the copy of $S^2\times S^4$ we consider, and $f_{[i]}$ is defined as above). Moreover, we perturb it in order to get only one critical point of index $6$, instead of two critical points of index $6$ and one of index $5$. We make the same perturbation for the minimum.

Hence its critical points are $p_6$, $p_4^1$, $p_4^2$, $p_2^1$, $p_2^2$, $p_0$ where as before $i_f (p_i)=i$ and $p_i^j \in (S^2\times S^4)_{[j]}$. Anew, if $\alpha_i$ and $\beta_i$ denote the generators of $\mS_1(\Omega S^2)$ and $\mS_3(\Omega S^4)$ in $(S^2\times S^4)_{[i]}$ for $i=1$ and $2$, we obtain in particular
\begin{align*}
\del p_6 = (\alpha_1 \otimes p_4^1 + \alpha_2 \otimes p_4^2) + (\beta_1 \otimes p_2^1 + \beta_2 \otimes p_2^2) + \gamma \otimes p_0.
\end{align*}

The element $\alpha_i \otimes p_4^i$ belongs to the second page of the Morse spectral sequence and, in view of (\ref{c2}), $\overline{c}^2(\alpha_1\otimes p_4^1)= f_{[1]}(p_4^1)$. Moreover $\del^2 (\alpha_i\otimes p_4^i)=0$ and $\del^2 p_6= \alpha_1 \otimes p_4^1 + \alpha_2 \otimes p_4^2$, hence at page $3$, $[\alpha_1 \otimes p_4^1]=[\alpha_2 \otimes p_4^2]\neq 0$ and then
\begin{align}
\overline{c}^3(\alpha_1\otimes p_4^1)= \min\{ f_{[1]}(p_4^1), f_{[2]}(p_4^2) \}.
\end{align}
It is therefore clear that we can now easily produce an example with $\overline{c}^3(\alpha_1\otimes p_4^1)<\overline{c}^2(\alpha_1\otimes p_4^1)$ by perturbing $f_{[2]}$ near $p_4^2$.

\begin{rema}
This leads obviously to an example in the Lagrangian intersection setting, by using the isomorphism due to Floer \cite{Floer89a}. For a manifold $L$ endowed with a Morse--Smale pair $(f,g)$ ($f$ is also required to be $C^2$--small), this isomorphism is induced by an identification between the two complexes
\begin{align*}
CM^*(L;f,g) \simeq CF^*(L,\Gamma_{df};f\circ\pi,J^g)
\end{align*}
with $\Gamma_{df}$ the graph of $df$ and $L$ seen as the $0$--section, Lagrangian submanifolds in $T^*L$. The map $\pi$ is the projection of the cotangent bundle and $J^g$ an almost complex structure induced by $g$. This isomorphism identifies the Morse and Floer spectral sequences of Barraud and Cornea. Moreover the action functional corresponds to the Morse function (up to an additive constant).
\end{rema}

\subsection{Distinguishing higher order spectral invariants}

In this section, we prove Theorem \ref{theo:main}. Recall from the introduction that this theorem provides a way to distinguish Lagrangian spectral invariants of higher order one from the other as soon as there exists a non trivial differential in $E_X(L)$. The difference between two such invariants is bounded in terms of the geometric quantity $r(L,L')$ which was defined in the introduction (Definition \ref{def:rll}).

\begin{rema}
As we have a large choice of spaces $X$, this theorem leads to many different examples. In particular if we could put $X=L$ (and thus, if we assume that $L$ is simply-connected), the second page of $E_L(L)$ would be non trivial. Moreover this spectral sequence would converge to the homology of a contractible space. Therefore there would be many non vanishing differentials.
\end{rema}

\begin{proof}[Proof of Theorem \ref{theo:main}]
Let $H$ be a Hamiltonian function such that $L=\phi^1_H(L')$ and $J$ be any almost complex structure such that the pair $(H,J)$ is regular. Recall that we denoted by $\Phi_f^H$ the extension of the PSS morphism to the Barraud--Cornea spectral sequences. We put $a=\Phi_f^H(\alpha)\in EF^r_{p,q}(L,L;H,J;X)$ and $b=\Phi_f^H(\beta)\in EF^r_{p-r,q+r-1}(L,L;H,J;X)$. They satisfy $d^r a=b\neq 0$.

Let $\eps$ be a small positive real number. We put $\nu=\skc{r}{\beta}-\eps$. By definition of $\skc{r}{\beta}$, we know that $q_\nu(b)\neq 0$. Since $q_\nu$ is a morphism of spectral sequences, $d^r (q_\nu (a))=q_\nu (d^r a)=q_\nu (b)$ and thus $q_\nu(a)\neq 0$. This proves the first inequality of the theorem for the relative version and thus for the spectral invariants.\\

Recall that the elements of $EF^r_{p,q}(L,L;H,J;X)$ consist of equivalence classes of elements of $F^p\mC^X_{p+q}$ with boundaries lying in $F^{p-r}\mC^X_{p+q-1}$. Hence, by definition of $\sic{r}{\alpha}$, we can choose such an element representing $a$, $\sum_{k\leq p} r_k\otimes \gamma_k$, where $r_k\otimes \gamma_k\in \mR_{p+q-k}\otimes CF_k$ and for every $k$, $\mA(\gamma_k)\leq \sic{r}{\alpha}+\varepsilon$. Assume that for every $\gamma\in\mI(L,L;\eta,H)$ such that $\mA(\gamma) \geq \skc{r}{\beta}-2\varepsilon$, $\mM_{\gamma_k,\gamma}(L,L;H,J)=\emptyset$ for all $k$. Then we can deduce that $q_\nu(a)$ is an $(r-1)$--cycle in the spectral sequence $EF^{>\nu}(L,L;H,J;X)$. Moreover its image by the $r$--th differential is well-defined and has to be $0$. This contradicts the fact that $q_\nu(b)\neq 0$.

Hence, there exist an element $\gamma_0\in\mI(L,L;\eta,H)$ and an index $k$ such that $\mA(\gamma_0) \geq \skc{r}{\beta}-2\varepsilon$ and $\mM_{\gamma_k,\gamma_0}(L,L;H,J)\neq\emptyset$. We get:
\begin{align*}
E(u) = \mA(\gamma_k)-\mA(\gamma_0) \leq \sic{r}{\alpha}-\skc{r}{\beta}+3\varepsilon.
\end{align*}
Introducing $\sigma_L(1;H,J,\eta)$ twice with opposite signs, we get
\begin{align*}
E(u) \leq \ic{r}{\alpha}-\kc{r}{\beta}
\end{align*}
since $\eps$ can be chosen arbitrarily small.

We recall that the naturality morphism maps a strip on a semi-tube with the same energy. Thus we get for all regular almost complex structure $J$, that there exist $x_J$ and $y_J\in L\cap L'$ and $u_J\in \mM_{x_J,y_J} (L,L';0,J)$ such that
\begin{align}\label{ineg1}
E(u_J)\leq \ic{r}{\alpha}-\kc{r}{\beta}.
\end{align}

Let us fix $\delta>0$ and denote $r_\delta:= r(L,L') -\delta$. We choose $\{ e^x_{r_\delta} |\, x\in L\cap L'\}$, a finite family of embeddings of the ball of radius $B(0,r_\delta)$, as in Definition \ref{def:rll} that is, satisfying (\ref{boules}):
\begin{align*}
{\rm i.} \;\; & (e^x_{r_\delta})^*(\omega) = \omega_0 \;\;\mbox{ and }\;\; e^x_{r_\delta}(0)=x, \\
{\rm ii.} \;\; & (e^x_{r_\delta})^{-1}(L) = \bR^n\cap B(0,{r_\delta})\;\; \mbox{ and } \;\; (e^x_{r_\delta})^{-1}(L') = i\bR^n\cap B(0,{r_\delta}).
\end{align*}
We let $J^\delta$ be an almost complex structure which coincides on ${\rm Im} \, e^x_{r_\delta}$ with $(e^x_{r_\delta})_*J_0$. If $J^\delta$ is not regular, we can find a sequence of regular almost complex structures, $\{ J_n \}_n$, converging to $J^\delta$. The above process leads to a sequence of pairs of orbits $\{ (x_n,y_n) \}_n$ and a sequence of Floer trajectories $\{ u_n \}_n$ (for all $n$, $u_n \in \mM_{x_n,y_n} (L,L';0,J_n)$). As $L\cap L'$ is a finite set, $\{ (x_n,y_n) \}_n$ admits a constant subsequence. We restrict ourselves to this subsequence $(x,y)$. The corresponding subsequence of $\{ u_n \}_n$, is a sequence of Floer trajectories $\{ u_k \}_k\subset \mM_{x,y} (L,L';0,J_k)$ with the same ends $x$ and $y$. Since the energy of Floer trajectories only depends on their ends, all the $u_k$'s have the same energy. Hence Gromov's compactness Theorem implies the existence of a $J^\delta$--pseudo-holomorphic limit strip $u_\delta$ whose energy satisfies (\ref{ineg1}).

By definition of the $e^x_{r_\delta}$'s, they preserve the symplectic area and thus
\begin{align*}
{\rm Area}_\omega({\rm Im}\, u_\delta \cap {\rm Im}\, e) = {\rm Area}_{\omega_0}( e^{-1} ({\rm Im}\, u_\delta \cap {\rm Im}\, e))
\end{align*}
where $e$ may be either $e^x_{r_\delta}$ or $e^y_{r_\delta}$. Moreover, by the choice of $J^\delta$, $e^{-1} ({\rm Im}\, u_\delta \cap {\rm Im}\, e)$ is a $J_0$--pseudo-holomorphic curve whose boundary lies in $\bR^n \cup i\bR^n \cup \del B(0,r_\delta)$. This extends, by symmetries to a pseudo-holomorphic curve, containing 0 in its interior, whose boundary lies in $\del B(0,r_\delta)$. The initial area has been multiplied by $4$. The isoperimetric inequality implies that the extended curve has area at least $\pi r_\delta^2$. Hence we get, ${\rm Area}_\omega({\rm Im}\, u_\delta \cap {\rm Im}\, e) \geq \frac{\pi}{4}\,r_\delta^2$.

Since this is true for both ends and ${\rm Im}\, e^x_{r_\delta} \cap {\rm Im}\, e^y_{r_\delta} = \emptyset$, we have
\begin{align*}
E(u_\delta)\geq {\rm Area}_\omega({\rm Im}\, u_\delta \cap {\rm Im}\, e^x_{r_\delta}) + {\rm Area}_\omega({\rm Im}\, u_\delta \cap {\rm Im}\, e^y_{r_\delta}) \geq \frac{\pi}{2}\,r_\delta^2.
\end{align*}
Therefore, for every (small) $\delta>0$,
\begin{align*}
\ic{r}{\alpha}-\kc{r}{\beta}\geq \frac{\pi}{2}\,(r(L,L') -\delta)^2.
\end{align*}
This ends the proof of Theorem \ref{theo:main}.
\end{proof}


\section{Back to homological Lagrangian spectral invariants}\label{sec:4}

Theorem \ref{theo:main} has a few interesting corollaries relative to the invariants of Definition \ref{defi:si}. We prove two of them in Section \ref{sec:coros} (Corollaries \ref{theo:ss} and \ref{prop:prop}) and then we show Proposition \ref{prop:dual} in Section \ref{sec:prop3}. In Section \ref{sec:fin}, we illustrate some other potential applications of the homological Lagrangian spectral invariants.

\subsection{Proof of Corollaries \ref{theo:ss} and \ref{prop:prop}}\label{sec:coros}

Corollary \ref{theo:ss} is just the homological counterpart of Theorem \ref{theo:main}.

\begin{proof}[Proof of Corollary \ref{theo:ss}]
If we assume that $X$ is $(r-1)$--connected, the pages $s$ with $2\leq s\leq r$ are identical to page $2$. Hence $E^r_{p,q}(L;X)\simeq H_q(\Omega X)\otimes H_p(L)$. We assume that $d^r \alpha =\beta\neq 0$ where $\alpha\in H_p(L)$ is identified with $1\otimes\alpha\in E^r_{p,0}(L;X)$ and $\beta$ denotes $\sum_i x_i\otimes \beta_i$ with $\beta_i$ in $H_{p-r}(L)$. Notice that once the ordered basis of $H_{r-1}(\Omega X)$, $\{x_i\}$, is chosen, the $\beta_i$'s are uniquely determined. Hence, we get as an immediate corollary that
\begin{align*}
\forall i,\; c(\alpha;L,L')-c(\beta_i;L,L') \geq \frac{\pi r(L,L')^2}{2}
\end{align*}
since $\ic{r}{1\otimes\alpha}=c(\alpha;L,L')$ and $\kc{r}{\beta}=\max_i \{ c(\beta_i;L,L') \}$. Corollary \ref{theo:ss} is proved.
\end{proof}

We now give the proof of Corollary \ref{prop:prop}. For $\alpha\in H_k(L)$, we denote by $\overline{\alpha}\in H^k(L)$ its $\operatorname{Hom}$--dual and by $\tilde{\alpha}\in H^{n-k}(L)$ the Poincar\'e dual of $\alpha$. As before, $\alpha\in H_k(L)$ is identified with $1\otimes\alpha\in E^r_{k,0}(L;X)$.
\begin{proof}[Proof of Corollary \ref{prop:prop}]
Recall that for any abelian group $\pi$ and any integer $n$, one can construct the Eilenberg--MacLane space $K(\pi,n)$, such that all its homotopy groups are trivial except the $n$--th which is $\pi$. They satisfy $\Omega K(\pi,n) = K(\pi,n-1)$. Moreover, we have the isomorphism
\begin{align*}
H^n (X,\pi)\simeq [X;K(\pi,n)].
\end{align*}
For $\alpha\in H_{k}(L,\bZ_2)$, with $k>1$, its $\operatorname{Hom}$--dual $\overline{\alpha}\in H^k(L)$ corresponds to the homotopy class of a map $l\co L\ra K(\bZ_2,k)$. This map is such that $l_*(\alpha)=\iota_k$, the generator of $H_k(K(\bZ_2,k),\bZ_2)\simeq \bZ_2$. Since $H_{k-1}(\Omega K(\bZ_2,k),\bZ_2)\simeq H_k(K(\bZ_2,k),\bZ_2)$, there has to be in $E_{K(\bZ_2,k)}(L)$ a non trivial differential
\begin{align*}
d^k \alpha=\iota_k\otimes 1\in H_{k-1}(\Omega K(\bZ_2,k), \bZ_2)\otimes H_0(L).
\end{align*}
Thus we conclude, via Corollary \ref{theo:ss}, that
\begin{align}
\frac{\pi r(L,L')^2}{2} \leq c(\alpha; L,L') - c(1; L,L') = c(\alpha; L,L')\label{preuve}.
\end{align}

\hop There is a cohomological version of the Serre spectral sequence $E_X(L)$ such that
\begin{itemize}
\item[i.] the differential of its $r$--th page has bidegree $(r,1-r)$,
\item[ii.] its second page comes with a product which coincides (as we use $\bZ_2$ coefficients) with the cup-product on cohomology classes.
\end{itemize}
We use the same process as above, replacing $\alpha$ by the ${\rm Hom}$--dual of its Poincar\'e dual $\overline{\tilde{\alpha}}$ and $K(\bZ_2,k)$ by $K(\bZ_2,n-k)$. Hence in the cohomological version of $E_X(L)$, there is a non trivial differential $d^{n-k} \overline{\iota_{n-k}}\otimes 1=\tilde{\alpha}$, as soon as $k<n-1$. By ii., we know that for any cohomology class $\gamma$, $d^{n-k} (\overline{\iota_{n-k}}\otimes \gamma)=\tilde{\alpha}\cup\gamma$. Putting $\gamma=\overline{\alpha\cdot\beta}$ gives that
\begin{align*}
d^{n-k} (\overline{\iota_{n-k}}\otimes \overline{\alpha\cdot\beta})=\overline{\beta}.
\end{align*}
Hence, in $E_{K(\bZ_2,n-k)}(L)$, the differential satisfies $d^{n-k} \beta=\iota_{n-k}\otimes (\alpha\cdot\beta)$. Corollary \ref{theo:ss} implies that
\begin{align*}
c(\alpha\cdot\beta; L,L') \leq c(\beta; L,L') - \frac{\pi r(L,L')^2}{2}.
\end{align*}
Thus the first assertion of Corollary \ref{prop:prop} is proved. Putting $\beta=[L]$ in this inequality, together with (\ref{preuve}), give the second one. Corollary \ref{prop:prop} is proved.
\end{proof}

\begin{rema}
In the next section we sketch another proof of this corollary. As mentioned in the introduction this alternate proof extends the result in the sense that it remains true for any non zero homology class $\alpha\in H_k(L)$, with $0<k<n$. However, the cases $k=1$ and $k=n-1$ should be obtained with the previous method by using a system of local coefficients.
\end{rema}

\subsection{Proof of Proposition \ref{prop:dual}}\label{sec:prop3}

In this section, we prove the properties of the homological Lagrangian spectral invariants listed in Proposition \ref{prop:dual}. Then we sketch an alternate proof of Corollary \ref{prop:prop}, via the moduli spaces used in the definition of the module structure.

\begin{proof}[Proof of Proposition \ref{prop:dual}]
The first assertion follows from the fact that there are Morse and Lagrangian forms of duality which are compatible with the PSS morphism. Indeed, given a Hamiltonian function $H$, on $M$, we denote by $H'$ the function defined by $H'(t,x):=-H(1-t,x)$. We have an identification
\begin{align*}
CF_*(L,L;H,J) \ni x \leftrightarrow x' \in CF_{n-*}(L,L;H',J)
\end{align*}
where $x$ and $x'$ are, geometrically, the same orbit with opposite orientation. This leads to a canonical isomorphism $HF_*(L,L;H,J) \simeq HF_{n-*}(L,L;H',J)$ and thus, via naturality,
\begin{align*}
HF_*(L,L';0,\widetilde{J}) \simeq HF_{n-*}(L,L'';0,\widetilde{J}')
\end{align*}
with $L'$ and $L''$ such that $\phi_H^1(L')=L$ and $\phi_H^1(L)=L''$. Moreover, this isomorphism preserves the absolute value of the action and changes its sign, that is for every $x\in L\cap L'$ and its corresponding $x'$, $\mA_{H'}(x')=-\mA_H(x)$ (this implies to choose corresponding references for the action functionals).

Given a Morse function $f$ on $L$, we denote $-f$ by $f'$. From the point of view of the moduli spaces defining the Lagrangian PSS morphism, we obviously get that $\mM_{p,x}^{f,H} = \mM^{x',p'}_{H',f'}$. Hence in homology $\phi_f^H(\alpha)=a$ if and only if $\phi_{f'}^{H'}(\alpha')=a'$ and $c(\alpha;L,L')=-c(\alpha';L,L'')$.

The second assertion comes from the two estimates:
\begin{align*}
\begin{split}
p\in {\rm Crit}_k f,\; \phi_f^{H}(p)=\sum_{\mu(y)=k} a_y y \mbox{ with } a_x\neq 0 &\Rightarrow \mA_H(x)\leq E_+(H)\\
y\in \mI(L,L;\eta,H),\; \psi_{H}^f(y)\neq 0 &\Rightarrow \mA_H(y)\geq E_-(H)
\end{split}
\end{align*}
relying on a computation (see~\cite{BarraudCornea05}) and leading to the following inequalities:
\begin{align*}
E_-(H)\leq \sigma_L(\alpha;H,J,\eta) + \alpha_\eta \leq E_+(H)
\end{align*}
for any $H$ such that $\phi^1_H(L)=L'$, $\alpha_\eta$ being a constant depending on our choice of reference $\eta$. Writing these inequalities for $\alpha=1$ the generator of $H_0(L)$ and subtracting them lead to
\begin{align*}
|c(\alpha;L,L')| \leq E_+(H)-E_-(H) = \| H \|.
\end{align*}
Since $c(\alpha;L,L')$ is non negative (this is a corollary of assertion {\it (2.)} of Corollary \ref{prop:prop}) and since the left side does not depend on $H$ as long as $\phi^1_H(L)=L'$, we get the desired inequality.
\end{proof}

Finally we sketch an alternate proof of Corollary \ref{prop:prop}, inspired by the Hamiltonian case \cite{Schwarz00} and using the definition of the module structure. Recall from the introduction that this proof extends Corollary \ref{prop:prop} in the sense that the condition on the degree of $\alpha$ is $0<k<n$ rather than $1<k<n-1$.

Indeed, the moduli spaces defining the module structure on $HF_*(L,L')$ are such that if $p$ is not the maximum of $f$, $\mM_{(p,x);y}\neq \emptyset$ implies that there exists a Floer trajectory connecting $x$ and $y$. The first step is to show that, in this case
\begin{align*}
\mA_H(x) \geq \mA_H(y) + \frac{\pi r(L,L')^2}{2}
\end{align*}
(one has to do the same procedure as at the end of the proof of Theorem \ref{theo:main}, with $J^\delta$). The second step consists in proving that
\begin{align*}
c(\alpha\cdot\beta; L,L') \leq c(\beta; L,L') - \frac{\pi r(L,L')^2}{2}
\end{align*}
as soon as $\alpha\cdot\beta\neq 0$ and $\alpha\neq [L]$ (this comes from the very definition of the spectral invariants and the previous inequality). Putting $\beta=[L]$, one gets
\begin{align*}
c(\alpha; L,L') \leq c([L]; L,L') - \frac{\pi r(L,L')^2}{2}
\end{align*}
for every non zero class $\alpha\in H_k(L)$, with $k<n$. With the notation of Proposition \ref{prop:dual}, we get for $\alpha\neq [L]'=1$ that
\begin{align*}
c(\alpha'; L,\phi^{-1}(L)) \leq c([L]; L,\phi^{-1}(L)) - \frac{\pi r(L,\phi^{-1}(L))^2}{2}
\end{align*}
where $\phi$ is a Hamiltonian diffeomorphism such that $\phi(L)=L'$. Then, the first assertion of Proposition \ref{prop:dual} gives that
\begin{align*}
c([L]; L,L') - c(\alpha; L,L') \leq c([L]; L,\phi^{-1}(L)) - \frac{\pi r(L,\phi^{-1}(L))^2}{2}
\end{align*}
and also implies in general that $c([L]; L,\phi^{-1}(L))=c([L]; L,\phi(L))$. Finally, since $r(L,\phi^{-1}(L))=r(L,\phi(L))$, we get the second inequality.

\subsection{Proof of Corollary \ref{coro:cup-length} and final remarks}\label{sec:fin}

Corollary \ref{coro:cup-length} follows immediately from Proposition \ref{prop:dual} and (the previous extension of) Corollary \ref{prop:prop}. It gives information on the geometry of the set of the Lagrangian submanifolds of $(M,\omega)$. Recall the definition of (the homological version of) the cup-length of $L$:
\begin{align*}
{\rm cl}(L):=\max\left\{ k+1 |\, \exists\, \alpha_i\in H_{d_i}(L),\, 1\leq i\leq k,\;\mbox{such that}\; \left.
\begin{array}{l}
0<d_i<n, \;{\rm and}\\
\alpha_1\cdot \ldots\cdot \alpha_k \neq 0
\end{array}
\right. \right\}.
\end{align*}
For such a maximal family $\{\alpha_i\}$, we have
\begin{align*}
c([L];L,L') = \sum_{1\leq l\leq {\rm cl}(L)} c(\beta_{l-1};L,L')-c(\beta_{l};L,L')
\end{align*}
where $\beta_0=[L]$, $\beta_l=\alpha_1\cdot \ldots\cdot \alpha_l$ ($1\leq l\leq {\rm cl}(L)-1$) and $\beta_{{\rm cl}(L)}=1$. Corollary \ref{prop:prop} gives then
\begin{align*}
c([L];L,L') \geq {\rm cl}(L)\cdot\frac{\pi r(L,L')^2}{2}.
\end{align*}
We conclude the proof of Corollary \ref{coro:cup-length} with the last assertion of Proposition \ref{prop:dual}.\\

\begin{rema}
We recall that the Hamiltonian spectral invariants can be viewed as particular homological Lagrangian ones via the BPS isomorphism. The second assertion of Proposition \ref{prop:dual}, namely $|c([L];L,L')|\leq \nabla(L,L')$ provides in general a better upper bound than Hofer's norm for symplectomorphisms. Indeed, Ostrover~\cite{Ostrover03} pointed out the existence of a family $\{ \varphi_t \}$, $t\in [0,\infty)$ in $\Ham$ and a constant $c$ such that:
\begin{align*}
d(\operatorname{id} ,\varphi_t)\ra \infty \mbox{ as } t\ra\infty \;\;\;\;\mbox{ and }\;\;\;\; \nabla(\Gamma_{\operatorname{id}}, \Gamma_{\varphi_t})=c \mbox{ for all } t.
\end{align*}
For such a family, the difference between the Hamiltonian spectral invariants associated to $[M]$ and $1$ remains bounded, since
\begin{align*}
\rho([M];\varphi_t)-\rho(1;\varphi_t) = c([\Delta];\Delta,\Gamma_{\varphi_t}) \leq \nabla(\Delta ,\Gamma_{\varphi_t})=c.
\end{align*}
\end{rema}

\begin{rema}
Another immediate consequence of the assertion {\it (2.)} of Proposition \ref{prop:prop} is that as soon as $L$ and $L'$ are transverse, Hamiltonian isotopic Lagrangian submanifolds, $\nabla(L,L')>0$. This is easily seen to give another proof of the well-known fact that Hofer's distance for Lagrangian submanifolds is non degenerate.

Moreover, following Schwarz \cite{Schwarz00}, we can use the very same arguments in order to deduce some properties of the set of Lagrangian submanifolds Hamiltonian isotopic to a fixed one, in particular cases. For example, let $M$ be the $2$--dimensional torus $S^1\times S^1$ and $L_0$ be the Lagrangian submanifold $\{ 0 \}\times S^1$. We consider for any integer $k$, the autonomous Hamiltonian function $H^k$ on $M$ defined by $H_k(x,y)=k\sin (2\pi x)$. The upper bound given by the assertion {\it (2.)} of Proposition \ref{prop:prop} allows to conclude that, when $k$ converges to infinity, $\nabla(L,(\phi^1_{H^k})^{-1}(L))$ converges to infinity. Therefore, the diameter of the set of Lagrangian submanifolds Hamiltonian isotopic to $L$ in $M$, endowed with Hofer's distance, is infinite. 
\end{rema}

\end{document}